\documentclass[leqno, myheadings, twoside]{amsart}
\usepackage{amsfonts}
\usepackage{amsmath, amsthm, amssymb, amscd, amsxtra,graphicx}
\usepackage{latexsym, amsfonts}
\usepackage{url}
\usepackage{texdraw}
\usepackage{epsfig}
\def\1-1{(1.1)}

\makeatletter \@addtoreset{equation}{section}

\makeatletter \renewcommand{\@biblabel}[1]{#1.}

\theoremstyle{remark}

\begin{document}
\title [Sharp Weyl-type formulas] {Sharp Weyl-type formulas of the spectral functions for biharmonic
Steklov eigenvalues}
\author{Genqian Liu}


\maketitle Department of Mathematics, Beijing Institute of
Technology,
 Beijing 100081, the People's Republic of China.
 \ \
E-mail address:  liugqz@bit.edu.cn

\vskip 0.46 true cm

\vskip 15 true cm

\begin{abstract}   In this paper, by explicitly calculating the principal symbols of
   pseudodifferential operators and by applying H\"{o}mander's spectral function theorem,
    we obtain the Weyl-type asymptotic formulas with sharp remainder
    estimates for the counting functions of the two classes of
    biharmonic Steklov eigenvalues $\lambda_k$ and $\mu_k$ in
    a smooth bounded domain of a Riemannian manifold.
    This solves a longstanding challenging problem.
   \end{abstract}

\vskip 1.39 true cm

\section{ Introduction}

\vskip 0.45 true cm

Let $(\mathcal{M}, g)$ be a $C^\infty$ Riemannian manifold of
dimension $n$
 with a positive definite metric tensor $g$, and
let $\Omega\subset {\mathcal{M}}$ be a bounded domain with
$C^\infty$ boundary $\partial \Omega$. Assume $\varrho$ is a
non-negative bounded function defined on $\partial \Omega$.
   We consider the following two biharmonic Steklov eigenvalue problems:
        \begin{eqnarray} \label{1-1}  \left\{\begin{array}{ll}\triangle^2_g u=0
      \quad \;\; &\mbox{in}\;\; \Omega,\\
  u=0 \;\; & \mbox{on}\;\; \partial \Omega,\\
   \triangle_g u+\lambda \varrho \frac{\partial u}{\partial \nu}=0
   \;\;& \mbox{on}\;\; \partial \Omega \end{array} \right. \end{eqnarray}
        and  \begin{eqnarray} \label{1-2}  \left\{\begin{array}{ll}\triangle^2_g v=0
      \quad \;\; &\mbox{in}\;\; \Omega,\\
  \frac{\partial v}{\partial \nu}=0 \;\; & \mbox{on}\;\; \partial \Omega,\\
   \frac{\partial (\triangle_g v)}{\partial \nu}-\mu^3 \varrho^3 v=0
   \;\;& \mbox{on}\;\;  \partial \Omega, \end{array} \right. \end{eqnarray}
      where $\nu$ denotes the inward unit
   normal vector to $\partial \Omega$,  and $\triangle_g$ is the Laplace-Beltrami operator
  defined in local coordinates by the expression,
   \begin{eqnarray*} \triangle_g =\frac{1}{\sqrt{|g|}}\sum_{i,j=1}^n
   \frac{\partial}{\partial x_i} \left( \sqrt{|g|}\,g^{ij}
   \frac{\partial}{\partial x_j}\right).\end{eqnarray*}
Here $| g | : =  det(g_{ij})$ is the determinant of the metric
tensor, and $g^{ij}$ are the components of the inverse of the metric
tensor $g$.

\vskip 0.12 true cm

(1.1) and (1.2) are the biharmonic Steklov eigenvalue problems (see \cite{BFG}, \cite{FGW},
  \cite{Ku}, \cite{KS2}, \cite{Pa1}, \cite{TG} and \cite{WX}).
 In each of the two cases, the spectrum is discrete and we arrange the eigenvalues in
non-decreasing order (repeated according to multiplicity)
\begin{eqnarray*}  & 0<\lambda_1 \le \lambda_2 \le  \cdots \le
   \lambda_k\le\cdots,
   \end{eqnarray*}
\begin{eqnarray*}  & 0=\mu_1 \le \mu_2 \le  \cdots \le
   \mu_k\le\cdots.
   \end{eqnarray*}
The corresponding eigenfunctions on $\partial \Omega$ are expressed as
$\frac{\partial u_1}{\partial \nu}$,
     $\frac{\partial  u_2}{\partial \nu},
   \cdots, \frac{\partial u_k}{\partial \nu}, \cdots$;
 $v_1$, $v_2, \cdots$, $v_k, \cdots$.
It is clear that $\lambda_k$ and $\mu_k$ can be characterized variationally as
\begin{eqnarray*}
 \lambda_1=\frac{\int_\Omega |\triangle_g u_1|^2 dx} {
\int_{\partial \Omega} \varrho \big(\frac{\partial u_1}{\partial
\nu}\big)^2 ds} =
  \inf_{\underset{
  0\ne \frac{\partial w}{\partial \nu}\in L^2(\partial \Omega)}{w\in H_0^1 (\Omega)\cap H^2 (\Omega)
  }}\;\; \frac{\int_\Omega |\triangle_g w|^2 dx} {
\int_{\partial \Omega} \varrho \big(\frac{\partial w}{\partial \nu}\big)^2 ds},\;\;\qquad \qquad \\
 \lambda_k=\frac{\int_\Omega |\triangle_g u_k|^2 dx} {
\int_{\partial \Omega} \varrho \big(\frac{\partial u_k}{\partial
\nu}\big)^2 ds} =
  \max_{\underset{codim(\mathcal {F})=k-1}
  {\mathcal{F}\subset H_0^1 (\Omega)\cap H^2 (\Omega)}
}  \,\,\inf_{\underset{0\ne \frac{\partial w}{\partial\nu}\in
L^2(\partial \Omega)}{w\in \mathcal{F}}}\;\; \frac{\int_\Omega
|\triangle_g w|^2 dx} { \int_{\partial \Omega} \varrho
\big(\frac{\partial w}{\partial \nu}\big)^2 ds}, \quad \; k=2,3, 4,
\cdots
\end{eqnarray*}
and
\begin{eqnarray*}
 \mu_0=0,\quad \, \mu_1^3=\frac{\int_D |\triangle_g v_1|^2 dx} {
\int_{\partial D} \varrho^3 v_1^2 ds} =
  \inf_{\underset{\int_{\partial D} \varrho^3 w\, ds =0}
  {w\in H^2 (D)}} \,
  \frac{\int_D |\triangle_g w|^2 dx} {
\int_{\partial D} \varrho^3  w^2 ds},\;\;\qquad \qquad \\
 \mu_k^3=\frac{\int_D |\triangle_g v_k|^2 dx} {
\int_{\partial D} \varrho^3  v_k^2 ds} =
  \max_{\underset{\int_{\partial D}
 \varrho^3 w\,ds=0\},\; \;  codim(\mathcal {F})=k}
  {\mathcal{F}\subset \{w\big|w\in  H^2 (D),}}
   \inf_{w\in \mathcal{F}}
\frac{\int_D |\triangle_g w|^2 dx} { \int_{\partial D} \varrho^3 w^2
ds}, \quad  k=2,3,4, \cdots
\end{eqnarray*}
where
 $H^m(\Omega)$ is the Sobolev space, and where $dx$ and $ds$ are the Riemannian
elements of volume and area on $\Omega$ and $\partial \Omega$,
respectively.

 The boundary value problems (\ref{1-1})  and (\ref{1-2}) have very
interesting interpretations in theory of elasticity. We refer the
reader to \cite{FGW}, \cite{Pa1} and \cite{WX} for more details.  In
view of the important applications, one is interested in finding the
asymptotic formulas for $\lambda_k$ and $\mu_k$ as $k\to \infty$. Let us
introduce the counting functions $A(\tau)$ and $B(\tau)$ defined as the numbers
of eigenvalues $\lambda_k$ and $\mu_k$ less than or equal to a given $\tau$, respectively.
Then our asymptotic problems for the eigenvalues are reformulated as
the study of the asymptotic behavior of
 $A(\tau)$ and $B(\tau)$ as $\tau\to +\infty$.

 The simpler harmonic Steklov problem was first introduced by
 V. A. Steklov for bounded domains in the plane in \cite{St}.  This problem is to
 find function  $v$ satisfying
\begin{eqnarray} \label{1-3}  \left\{\begin{array}{ll} \triangle_g v=0&
\quad \, \mbox{in}\;\; \Omega,\\
\frac{\partial v}{\partial \nu}+ \eta \varrho v =0 &\quad \,
\mbox{on}\;\; \partial \Omega,
\end{array} \right.\end{eqnarray}
 where $\eta$ is a real number (The function $v$ represents
 the steady state temperature on $\Omega$ such that
 the flux on the boundary is proportional
 to the temperature). The harmonic Steklov spectrum of the domain is also called as the spectrum
of the Dirichlet-to-Neumann map (see \cite{Cal},
 \cite{FS}) or \cite{SU}.  For the harmonic Steklov eigenvalue problem (\ref{1-3}), in 1955
  Sandgren \cite{Sa}  established the asymptotic formula of
 the counting function $N(\tau)=
 \#\{k\big|\eta_k \le \tau\}$:
   \begin{eqnarray} \label{1-6}   N(\tau)=
   \frac{\omega_{n-1}\tau^{n-1}}{(2\pi)^{n-1}}
 \int_{\partial \Omega} \varrho^{n-1} ds +o(\tau^{n-1}) \quad \;\mbox{as}\;\;
 \tau\to +\infty,\end{eqnarray}
 where $\omega_{n-1}$ is the volume of the unit ball of ${\Bbb R}^{n-1}$.
  In the case that Riemannian manifold $\mathcal{M}$ and the boundary
 of $\Omega$ are smooth,
 the author \cite{Liu1} further gave a sharp remainder
 estimate for the counting function of
  the harmonic Steklov eigenvalues $\{\eta_k\}_{k=1}^\infty$:
\begin{eqnarray*}
N(\tau)=\frac{\omega_{n-1}\tau^{n-1}}{(2\pi)^{n-1}}
 \int_{\partial \Omega} \varrho^{n-1}(x) dx +O(\tau^{n-2})
 \quad \;
\mbox{as}\;\; \tau\to +\infty.\end{eqnarray*} For the biharmonic
 Steklov eigenvalue problem (\ref{1-1}) with general domain, in \cite{Liu1} the author also established the
leading asymptotic formula with remainder $o(\tau^{n-1})$ as
$\tau\to +\infty$.

 \vskip 0.16 true cm

For the Dirichlet and Neumann eigenvalues of the Laplacian on a bounded domain
$\Omega\subset {\Bbb R}^n$,
 H. Weyl (\cite{We1}, \cite{We2}) in 1912 proved the following
 asymptotic formula which answered a question posed in 1908 by the
 physicist Lorentz:
\begin{eqnarray} \label{1,,1} N_{\mp}(\tau)= (2\pi)^{-n} \omega_{n}
(\mbox{vol}(\Omega))\tau^{n/2}+o(\tau^{n/2})\quad \;\,
\mbox{as}\;\; \tau\to +\infty,\end{eqnarray} where
  $N_{-}(\tau)=\#\{k \in {\Bbb N}\big|\alpha_k\le \tau\}$,
  $N_{+}(\tau)=\#\{k \in {\Bbb N}\big|\beta_k\le \tau\}$,
  and $0<\alpha_1<\alpha_2 \le \cdots \le \alpha_k\le \cdots$
  and  $0=\beta_1<\beta_2 \le \cdots \le \beta_k\le \cdots$
 are all the Dirichlet and Neumann eigenvalues on $\Omega$, respectively.
  As far back as in 1912, H. Weyl \cite{We5} conjectured (see also,
Clark \cite{Cl})
 that the second term of the asymptotic formula for $N_{\mp}(\tau)$ contain
 an $(n-1)$-dimensional measure (`area')
of the boundary $\partial \Omega$, i.e.,
\begin{eqnarray} \label{1,,2} N_{\mp}(\tau)&=& (2\pi)^{-n}
\omega_{n} (\mbox{vol}(\Omega)) \tau^{n/2} \mp \frac{1}{4}
(2\pi)^{-n+1}
\omega_{n-1} (\mbox{vol}(\partial \Omega)) \tau^{(n-1)/2}\\
&& \quad \;\quad \; +o(\tau^{(n-1)/2}) \quad \;\mbox{as}\;\;
\tau\to +\infty.\quad \nonumber\end{eqnarray}
In 1980, Ivrii \cite{Iv} proved this conjecture for domains having smooth boundary
 under the following condition regarding the billiard trajectories of $\Omega$, where
 the billiard trajectory in $\Omega$ is taken with the
 usual reflections at the boundary (see also p.$\;$100 of \cite{ES}).
  Melrose \cite{Me} independently obtained the second term
 in Weyl's conjecture for manifolds with concave boundary.

 Note that it was already observed by Avakumovi\v{c} \cite {Av}
that for the Laplacian on the sphere ${\Bbb S}^n$, the high
multiplicities of the eigenvalues make it impossible to improve
(\ref{1,,1}) to (\ref{1,,2}).
Seeley (see, \cite{Se1} and \cite{Se2}) in 1980 gave the
sharp asymptotic formula of the counting function (see also \cite{ANPS}, \cite{ES}, \cite{SV}):
\begin{eqnarray} \label{1,,3} N_{\mp}(\tau)= (2\pi)^{-n} \omega_{n}
(\mbox{vol}(\Omega))\tau^{n/2}+O(\tau^{(n-1)/2})\quad \;\,
\mbox{as}\;\; \tau\to +\infty.\end{eqnarray}
 Applying the sharp
asymptotic result (\ref{1,,3}), Sogge invented the unit band spectral
projection operator (see \cite{So} and \cite{SS}) and
 established the well-known asymptotic estimates of the mapping norm
$\|\chi_\tau\|_{L_2\to L_p} \,\, (2\le p\le \infty)$ (cf.$\,$\cite{SZ}). Howerver,
it has been a longstanding challenging problem
to get the sharp Weyl-type asymptotic formulas for the biharmonic Steklov
 eigenvalues (see \cite{Liu1}).
\vskip 0.25 true cm

 In this paper, by explicitly calculating the principal symbols of
the corresponding pseudodifferential operators for the problems (\ref{1-1}) and (\ref{1-2})
in the boundary of
a $C^\infty$ bounded
domain, we obtain the sharp asymptotic formulas for the
counting functions $A(\tau)$ and $B(\tau)$, respectively. The main results are the
following:

\vskip 0.25 true cm

 \noindent  {\bf Theorem 1.1.} \   {\it
Let $(\mathcal{M},g)$ be an $n$-dimensional $C^\infty$ Riemannian
manifold, and let $\Omega\subset \mathcal{M}$ be a bounded domain
with $C^{\infty}$ boundary $\partial \Omega$. Then
\begin{eqnarray} \label{1-7} A(\tau)=
\frac{\omega_{n-1}\tau^{n-1}}{(4\pi)^{n-1}}
 \int_{\partial \Omega} \varrho^{n-1} ds +O(\tau^{n-2})\quad \,
 \; \mbox{as}\;\; \tau\to +\infty.\end{eqnarray}
 Moreover,  the above remainder
estimate is sharp}.

\vskip 0.25 true cm

 \noindent  {\bf Theorem 1.2.} \   {\it
Let $(\mathcal{M},g)$ be an $n$-dimensional $C^\infty$ Riemannian
manifold, and let $\Omega\subset \mathcal{M}$ be a bounded domain with
smooth boundary $\partial \Omega$.  Then
\begin{eqnarray} \label{1;-6;} B(\tau) = \frac{\omega_{n-1}\tau^{n-1}}{(\sqrt[3]{16}\pi)^{n-1}}
 \int_{\partial D} \varrho^{n-1} ds + O( \tau^{n-2})\quad \, \; \mbox{as}\;\; \tau\to +\infty.\end{eqnarray}
 Moreover,  the above remainder
estimate is also sharp}.

\vskip 0.25 true cm

 The plan of the paper is as follows. In Section
2 we give some definitions and lemmas. In Section 3, by a key technique we calculate
the principal symbols of the corresponding ``Neumann-to-Laplacian
map'' and ``Dirichlet-to-Laplacian derivative map''.
 Section 4 is devoted to
 the proofs of the sharp Weyl-type
 asymptotic formulas for $A(\tau)$ and $B(\tau)$. In
Section 5, we give two counterexamples, which show that Theorem
1.1 and 1.2 cannot be improved.

\vskip 1.39 true cm

\section{Some definitions and lemmas}

\vskip 0.45 true cm

 \noindent  {\bf Definition 2.1.} \ \  {\it  If
$U$ is an open subset  of ${\Bbb R}^n$, we denote  by $S^m=S^m (U,
{\Bbb R}^n)$ the set of all $p\in C^\infty (U, {\Bbb R}^n)$ such
that for every compact set $K\subset U$
 we have
 \begin{eqnarray} \label {2.3} |D^\beta_x D^\alpha_\xi p(x,\xi)|\le C_{K,\alpha,
 \beta}(1+|\xi|)^{m-|\alpha|}, \quad \; x\in K,\,\, \xi\in {\Bbb R}^n\end{eqnarray}
 for all $\alpha, \beta\in {\Bbb N}^n_+$.
 The  elements  of $S^m$  are  called  symbols of order $m$}.

\vskip 0.13 true cm

 It is clear that $S^m$ is a
Fr\'{e}chet space with semi-norms given by the smallest constants
which can be used in (\ref{2.3}) (i.e.,
\begin{eqnarray} \|p\|_{K,\alpha, \beta}=
 \,\sup_{x\in K}\bigg|\left(D_x^\beta D_\xi^\alpha
 p(x, \xi)\right)(1+|\xi|)^{|\alpha|-m}\bigg|).\end{eqnarray}

 Let $p(x, \xi)\in S^m$. A pseudo-differential operator
in an open set $U\subset {\Bbb R}^n$ is essentially defined by a
Fourier integral  operator (cf. \cite{Ho4}):
\begin{eqnarray} \label{2.1}  P(x,D) u(x) = \frac{1}{(2\pi)^{n}} \int_{{\Bbb R}^n} p(x,\xi)
 e^{i\langle x,\xi\rangle} \hat{u} (\xi)d\xi.\end{eqnarray}
Here $u\in C_0^\infty (U)$ and $\hat{u} (\xi)=
 \int_{{\Bbb R}^n} e^{-i\langle y, \xi\rangle} u(y)dy$ is the Fourier
transform of $u$.

\vskip 0.29 true cm

 \noindent  {\bf Definition 2.2.} \ \  {\it
 A pseudodifferential
operator $P$ with its symbol $p$ in $S^m$ is called {\it classical
or polyhomogeneous} if there is a sequence of symbols $p_j\in
S^{m-j}$, $\,j=0,1,2,\cdots,$ such that $p_j(x, t\xi)=t^{m-j}p_j (x,
\xi)$ for $t>1$, $|\xi|>1$, and
\begin{eqnarray*} \bigg|D^\beta_x D^\alpha_\xi \bigg(p(x, \xi)- \sum_{j=0}^N p_j(x, \xi)\bigg)\bigg|
\le C_{\alpha, \beta, N} |\xi|^{m-N-1-|\alpha|}\end{eqnarray*} for
all $\alpha, \beta, |\xi|>1$ and all integers $N\ge 0$. In this case
the notation $p\sim\sum_{j=0}^\infty p_j$ is used. The function
$p_0$ is known as the {\it principal symbol} of pseudodifferential
 operator $P$, and the class of such symbol is denoted by
$S^m_{cl}$.}

\vskip 0.19 true cm

Given a diffeomorphism
  $\iota:  U_1 \to U_2$, from one open set $U_1\subset {\Bbb R}^n$ onto another open
set $U_2\subset {\Bbb R}^n$, the induced transformation
$\iota^* : C^\infty_0 (U_2)\to C^\infty_0 (U_1)$,  taking a
function $u$ to the function $u\circ \iota$, is an isomorphism and
transforms $C_0^\infty(U_2)$ into $C_0^\infty(U_1)$. Let $P_1$ be
a pseudodifferential operator on $U_1$ and define $P_2:
C_0^\infty (U_2)\to C_\infty(U_2)$
with the help of the commutative diagram
\begin{eqnarray*}
\begin{CD}
C_0^\infty (U_1) @> P_1>>  C^\infty(U_1) \\
@A\iota^* AA @AA\iota^* A \\
 C_0^\infty(U_2) @> P_2>> C^\infty(U_2)
\end{CD}
\end{eqnarray*}
 i.e.,
  \begin{eqnarray} \label {2.12} P_2 u  = [P_1 (u \circ \iota)] \circ \iota^{-1}.\end{eqnarray}
(\ref{2.12}) can also be written as
 \begin{eqnarray*}  P_2 u  = (\iota^{-1})^* P_1 (\iota^* u).\end{eqnarray*}
It follows from this that $P_2$ is also a pseudodifferential
operator on $U_2$.

\vskip 0.1 true cm

 Let $\mathcal{M}$ be a smooth $n$-dimensional Riemannian
 manifold (of class $C^\infty$). We will denote by $C^\infty(\mathcal{M})$ and $C_0^\infty(\mathcal{M})$
 the space of all smooth complex-valued functions on $\mathcal{M}$  and the subspace of all
 functions with compact support, respectively. Assume that we are given a linear
 operator
 \begin{eqnarray*} P: C^\infty_0(\mathcal{M}) \to C^\infty(\mathcal{M}).\end{eqnarray*}
 If $G$ is some chart in $\mathcal{M}$ (not necessarily connected) and $\kappa: G\to U$ its
diffeomorphism onto an open set $U \subset {\Bbb R}^n$, then let ${\tilde P}$ be defined by the diagram
\begin{eqnarray*}
\begin{CD}
C_0^\infty (G) @> P>>  C^\infty(G) \\
@A\kappa^* AA @AA\kappa^* A \\
 C_0^\infty(U) @> \tilde P>> C^\infty(U)
\end{CD}
\end{eqnarray*}
  (note, in the upper row is the operator $r_G\circ P\circ i_G$, where $i_G$ is the natural
embedding $i_G: C_0^\infty (G)\to C^\infty_0 (M)$ and $r_G$ is the natural restriction $r_G: \, C^\infty(M)\to
C^\infty(G)$; for brevity we denote this operator by the same letter $P$ as the
original operator).

\vskip 0.29 true cm

\noindent  {\bf Definition 2.3.} \  {\it
 An operator $P: C_0^\infty (\mathcal{M}) \to C^\infty(\mathcal{M})$ is called a pseudodifferential
operator on $\mathcal{M}$ if for any chart diffeomorphism $\kappa: G\to U$, the
operator $\tilde P$ defined above is a pseudodifferential
operator on $U$}.

\vskip 0.29 true cm

 \noindent  {\bf Lemma 2.4 (see, for example,
 Proposition 0.3.C of \cite{Ta3})} \ \  {\it  If $A$ and $B$ are two pseudodifferential operators
of order $m$ and $m'$, respectively, then the composition $C=A\circ B$
is a pseudodifferential operator of order $m+m'$
 with the symbol
\begin{eqnarray} c(x, \xi) \sim \sum_{\alpha} \frac{i^{|\alpha|}}{\alpha!} D^\alpha_\xi
a(x, \xi)  D^\alpha_x b(x, \xi)\end{eqnarray} where $a(x, \xi)$ and
$b(x, \xi)$ are the symbols of $A$ and $B$, respectively. In
particular, the principal symbol of $A\circ B$ is
$a_0(x,\xi)b_0(x,\xi)$, where $a_0(x, \xi)$ and $b_0(x, \xi)$  are
the principal symbols of $A$ and $B$, respectively.}

\vskip 0.26 true cm

\noindent {\bf Lemma 2.5.} \ \ {\it  Let
  \begin{gather} \label{2--1}  A=\begin{pmatrix} a^{11} & a^{12} & \cdots  & a^{1,n-1} &0\\
      a^{21} &  a^{22} & \cdots & a^{2, n-1} & 0\\
      \vdots & \vdots & \ddots & \vdots & \vdots\\
      a^{n-1, 1} & a^{n-1,2} & \cdots & a^{n-1, n-1}  & 0\\
      0 & 0 & \cdots   & 0 &  a^{nn}  \end{pmatrix} \end{gather}
be a positive definite, real symmetric constant matrix.
 Let $\phi(x')$ and $h(x')$ be $C^\infty$ functions
of compact support in $(n-1)$-space.
Then the problem
\begin{eqnarray} \label {2-0-1} \left\{ \begin{array}{ll} \left(
\sum_{j,k=1}^{n-1} a^{jk} \frac{\partial^2 }{\partial x_j \partial
x_k}
 + a^{nn}\frac{\partial^2 }{\partial x_n^2}\right)^2 u =0 \;\quad \quad &\mbox{in}\;\; {\Bbb R}^n_+, \\
  u=\phi, \, \quad \quad \quad & \mbox{on}\;\; \partial {\Bbb R}^n_+, \\
 \sqrt{a^{nn}}\, \frac{\partial u}{\partial  x_n}=h \quad  \quad
 & \mbox{on}\;\; \partial {\Bbb R}^n_+\end{array}\right. \end{eqnarray}
 has a solution
\begin{eqnarray} \label {2-0-2} \quad \quad u(x',x_n) =\int_{{\Bbb R}^{n-1}} K_1(x'-y', x_n)
 \phi(y') dy' + \int_{{\Bbb R}^{n-1}}  K_2(x'-y', x_n)
h(y') dy',\end{eqnarray}
  where  $\,{\Bbb R}^n_+=\{x=(x_1, \cdots, x_{n-1},
x_n)\in {\Bbb R}^n \big|x_n>0\}$, and
 \begin{eqnarray*}   K_1(x', x_n) &=&
  (-1)^{n-1}\frac{(n-2)!}{(2\pi i)^{n-1}}\int_{|\eta'|=1}
 \left[\left(x'\cdot\eta' +i x_n \sqrt{\sum_{j,k=1}^{n-1}
 \frac{a^{jk}\eta_j\eta_k}{a^{nn}}}\right)^{1-n} \right.\\
 && \left. \,\, + (n-1)i x_n \sqrt{\sum_{j,k=1}^{n-1} \frac{a^{jk}\eta_j\eta_k }{a^{nn}}}
  \left(x'\cdot \eta' +i x_n
 \sqrt{\sum_{j,k=1}^{n-1}
 \frac{a^{jk}\eta_j\eta_k}{a^{nn}}}\right)^{-n}\right]ds_{\eta'},\\
   K_2(x', x_n) &=&
  (-1)^{n-1} \frac{(n-2)!}{(2\pi i)^{n-1}}\int_{|\eta'|=1}  \left[
 \frac{x_n}{\sqrt{a^{nn}}} \left(x'\cdot \eta'
 +i x_n \sqrt{\sum_{j,k=1}^{n-1} \frac{a^{jk}\eta_j \eta_k}{a^{nn}}}
 \right)^{1-n}\right] ds_{\eta'}.\end{eqnarray*}
 Here $\eta'=(\eta_1, \cdots, \eta_{n-1})$ and
$\,ds_{\eta'}$ is the area
 element on the unit sphere $|\eta'|=1$.}

\vskip 0.36 true cm

\noindent  {\it Proof.} \  Writing $x=(x',x_n)$.
 Then the bi-Laplace operator $P$ has characteristic form
  $P(\eta', \tau)=(\sum_{j,k=1}^{n-1} a^{jk}\eta_j\eta_k +a^{nn}\tau^2)^2$.
 It is easy to see that the roots of
   $P(\eta', \tau)$ with positive imaginary parts are $\tau^+_1 (\eta')= \tau^+_2(\eta')
   =i\,\sqrt{\sum_{j,k=1}^{n-1}\frac{a^{jk}\eta_j\eta_k}{a^{nn}}}$.
   Thus we have (see, Chapter I, \S1 of \cite{ADN})
\begin{eqnarray*} M^+(\eta', \tau) = \left(\tau-i\,\sqrt{\sum_{j,k=1}^{n-1}
\frac{a^{jk}\eta_j\eta_k}{a^{nn}}}\right)^2
 = \tau^2 -2i
\left(\sqrt{\sum_{j,k=1}^{n-1}\frac{a^{jk}\eta_j\eta_k}{a^{nn}}}\right)\tau
- \sum_{j,k=1}^{n-1}
\frac{a^{jk}\eta_j\eta_k}{a^{nn}},\end{eqnarray*}
\begin{eqnarray*} M^+_0(\eta', \tau) = 1,  \;\quad \quad  M^+_1(\eta', \tau) =
\tau -2i\,\sqrt{\sum_{j,k=1}^{n-1} \frac{a^{jk}
\eta_j\eta_k}{a^{nn}}},\end{eqnarray*} so that
\begin{eqnarray*} N_1 (\eta', \tau)= M_1^+ (\eta', \tau) = \tau -2i\,\sqrt{\sum_{j,k=1}^{n-1}
\frac{a^{jk} \eta_j\eta_k}{a^{nn}}},\; \quad\;\;
 N_2 (\eta', \tau)= \frac{1}{\sqrt{a^{nn}}}M_0^+ (\eta', \tau) =\frac{1}{\sqrt{a^{nn}}}.\end{eqnarray*}
It follows from p.$\;$635 of \cite{ADN}  and the well-known
 residue theorem (see, for example, p.$\,$150 of \cite{Ah}) that
\begin{eqnarray*} & K_1(x', x_n) = (-1)^{n-1} \,\frac{(n-2)!}{(2\pi i)^n }
\int_{|\eta'|=1} \left[ \int_\gamma \frac{N_1(\eta',
\tau)}{
M^+(\eta', \tau)(x'\cdot \eta' +x_n\tau)^{n-1}} d\tau \right] ds_{\eta'}
\qquad\qquad\qquad \qquad\qquad \qquad \;\quad\\
 & = (-1)^{n-1} \,\frac{(n-2)!}{(2\pi i)^n}
\int_{|\eta'|=1} \left[ \int_\gamma \frac{\tau- 2i
\sqrt{\sum_{j,k=1}^{n-1} \frac{a^{jk} \eta_j\eta_k}{a^{nn}}}}
 {\left(\tau-i\,\sqrt{\sum_{j,k=1}^{n-1}\frac{a^{jk}
 \eta_j\eta_k}{a^{nn}}}\right)^2
 (x'\cdot \eta'+x_n \tau)^{n-1}} d\tau \right]ds_{\eta'}\quad\quad \\
 & = (-1)^{n-1}\frac{(n-2)!}{(2\pi i)^{n-1}}\int_{|\eta'|=1}
 \left[\left(x'\cdot\eta' +i x_n \sqrt{\sum_{j,k=1}^{n-1}
 \frac{a^{jk}\eta_j\eta_k}{a^{nn}}}\right)^{1-n} \right.\quad \quad\qquad \qquad \\
  &\;\; \left.  + (n-1)i x_n \sqrt{\sum_{j,k=1}^{n-1} \frac{a^{jk}\eta_j\eta_k }{a^{nn}}}
  \left(x'\cdot \eta' +i x_n
 \sqrt{\sum_{j,k=1}^{n-1}
 \frac{a^{jk}\eta_j\eta_k}{a^{nn}}}\right)^{-n}\right]ds_{\eta'},
 \quad  \end{eqnarray*}
 \begin{eqnarray*}
    & K_2(x', x_n) = (-1)^{n-2} \, \frac{(n-3)!}{(2\pi i)^n}
     \int_{|\eta'|=1}ds_{\eta'} \left[ \int_\gamma
\frac{N_2(\eta', \tau)}{ M^+(\eta', \tau)(x'\cdot \eta' +x_n\tau)^{n-2}}
d\tau \right]\qquad \qquad \qquad \;\qquad \qquad\quad\qquad\qquad\\
& =(-1)^{n-2} \; \frac{(n-3)!}{(2\pi i)^n}
 \int_{|\eta'|=1} \left[ \int_\gamma \frac{
1} {\sqrt{a^{nn}}\left(\tau-i \,\sqrt{\sum_{j,k=1}^{n-1}\frac{a^{jk}
\eta_j\eta_k}{a^{nn}}} \right)^2 (x'\cdot \eta'+x_n \tau)^{n-2}}
d\tau \right] ds_{\eta'}\quad \;\;\\
&= (-1)^{n-1} \frac{(n-2)!}{(2\pi i)^{n-1}}\int_{|\eta'|=1}  \left[
 \frac{x_n}{\sqrt{a^{nn}}} \left(x'\cdot \eta'
 +i x_n \sqrt{\sum_{j,k=1}^{n-1} \frac{a^{jk}\eta_j \eta_k}{a^{nn}}}
 \right)^{1-n}\right] ds_{\eta'},\qquad \;\; \;\end{eqnarray*}
 where $\gamma$ is a Jordan contour in $\mbox{Im} \,\tau>0$
 enclosing all the points $i\sqrt{\sum_{i,j=1}^{n-1} \frac{a^{jk}
 \eta_i\eta_j}{a^{nn}}}$ for all $|\eta'|=1$.
  Applying Theorem 2.1 of \cite{ADN}, we obtain (\ref{2-0-2}).  $ \quad \quad \square $

 \vskip 0.18 true cm

Let $\{E_\tau\}$ be the spectral resolution of pseudodifferential operator
 $P$, and let $e(x,y,\tau)$ be the kernel of $E_\tau$. This is an element of
 $C^\infty (\Omega\times \Omega)$  called the spectral function of $P$.

   \vskip 0.2 true cm

The following Lemma will be used later.

\vskip 0.2 true cm

 \noindent {\bf Lemma 2.6 (H\"{o}mander's spectral function theorem, see,
 Theorem 5.1 of  \cite{Ho}, \cite{Ho3} or \cite{Sh})} \ \  {\it
Let $P$ be a
 non-negative pseudodifferential operator, acting on a $C^\infty$ subdomain $\Omega$
  of an $n$-dimensional $C^\infty$ manifold. Let  $p_0(x,\xi)$  be
  the  principal  symbol of $P$,  which  is a real
homogeneous polynomial of degree $m$  on  the  cotangent  bundle
$T^*(\Omega)$.  The measure $dx$ defines  a  Lebesgue measure $d\xi^*$
in each fiber of  $T^*(\Omega)$;  which  is  a  vector space of
dimension $n$. Then \begin{eqnarray} \tau^{-n/m} e(x,x,\tau)
-(2\pi)^{-n} \int_{B_x} d\xi^* =O(\tau^{-1/m}) \quad
\,\mbox{as}\;\; \tau\to \infty,\end{eqnarray} where $B_x =\{
\xi\in T^*_x(\Omega) \big|p_0(x,\xi)<1\}$.
}

\vskip 1.39 true cm

\section{The principal symbols}

\vskip 0.45 true cm

3.1.  \ \  Let $(\mathcal{M}, g)$ be a $C^\infty$ Riemannian manifold, and let
$\Omega$ be a bounded domain with $C^\infty$ boundary in
$\mathcal{M}$.  The ``Neumann-to-Laplacian map'' is the map
\begin{eqnarray*} F: H^{1/2} (\partial \Omega)\to H^{-1/2} (\partial \Omega)\end{eqnarray*}
  defined by the following problem:  Let $h\in H^{1/2} (\partial \Omega)$ and
let $u\in H^2(\Omega)$ be the  solution of
\begin{eqnarray} \label{3.1}  \left\{ \begin{array} {ll}
\triangle_g^2 u=0 \,\quad\;\mbox{in}\;\; \Omega,\\
  u=0 \quad \;\; \quad \,\mbox{on}\;\; \partial \Omega, \\
 \frac{\partial u}{\partial \nu}=h \quad\quad \mbox{on}\;\; \partial
 \Omega,\end{array}\right. \end{eqnarray}
 we set $Fh:= (-\Delta_g u)\big|_{\partial \Omega}$.
    Multiplying (\ref{3.1}) by $u$, integrating the result over $\Omega$,
and using Green's formula, we derive
\begin{eqnarray*} 0&=&\int_\Omega u(\Delta_g^2 u)dx =\int_\Omega
|\Delta_g u|^2 dx -\int_{\partial \Omega}u\frac{\partial (\Delta_g
u)}{\partial \nu} ds +\int_{\partial \Omega} (\Delta_g
u)\frac{\partial u}{\partial \nu} ds\\
&=& \int_\Omega |\Delta_g u|^2 dx +\int_{\partial \Omega}
 (\Delta_g u)\frac{\partial u}{\partial \nu}ds,\end{eqnarray*} so that
  $$\langle Fh, h\rangle=\int_{\partial \Omega}
(Fh)h\, ds =\int_\Omega |\triangle_g u|^2 dx\ge 0,\;\; \,
\mbox{for any}\;\; h\in H^{1/2} (\partial \Omega). $$ This shows
that $F$ is a non-negative, self-adjoint, pseudodifferential
  operator on $H^{1/2}(\partial \Omega)$.
We shall  calculate the principal symbol of $F$.

\vskip 0.13 true cm

 \noindent {\bf Lemma 3.1} \ \ {\it  Let $A$ be a positive definition, real-valued constant
  matrix as in (\ref{2--1}).
 Assume that
\begin{eqnarray*} F_0: C^\infty_0({\Bbb R}^{n-1}) \to C^\infty ({\Bbb R}^{n-1})\end{eqnarray*}
  defined by the following problem:  Let $h\in C^\infty_0 ({\Bbb R}^{n-1})$ and
let $u\in C^\infty({\Bbb R}^n_+)$ be the  solution of
\begin{eqnarray} \label{3..2} \left\{\begin{array}{ll} \left( \sum_{j,k=1}^{n-1}
a^{jk} \frac{\partial^2}{\partial x_j \partial x_k} + a^{nn}
\frac{\partial^2 }{\partial x_n^2}\right)^2 u=0
\quad \quad \;\; & \mbox{in}\;\; {\Bbb R}^n_+, \\
u=0 \quad \quad \quad \quad\quad \quad & \mbox{on}\;\; \partial {\Bbb R}^n_+, \\
 \sqrt{a^{nn}}\frac{\partial u}{\partial x_n}=h \;\quad\;\quad
  & \mbox{on}\;\; \partial {\Bbb R}^n_+.\end{array}\right.\end{eqnarray}
 we set $F_0h:= - \left(\sum_{j,k=1}^{n-1}
a^{jk} \frac{\partial^2}{\partial x_j \partial x_k} + a^{nn}
\frac{\partial^2 }{\partial x_n^2}\right)\big|_{\partial {\Bbb R}^n_+}$.
   Then  the principal symbol of $F_0$ is
 \begin{eqnarray*}  p_0(x', \eta') =2
\sqrt{\sum_{j,k=1}^{n-1} a^{jk} \eta_j\eta_k},\quad \quad \forall
\,\,  (x', \eta')\in {\Bbb R}^{n-1}\times ({\Bbb R}^{n-1}\setminus 0).\end{eqnarray*}}

\vskip 0.2 true cm

\noindent  {\it Proof.} \ \   Writing  $x=(x', x_n)$,  it follows from Lemma 2.5 that
 \begin{eqnarray} \label {3..3} u(x', x_n)= \int_{{\Bbb R}^{n-1}} K_2(x'-y', x_n) h(y') dy',\end{eqnarray}
where
 \begin{eqnarray*}
K_2(x', x_n) =
  (-1)^{n-1} \frac{(n-2)!}{(2\pi i)^{n-1}}\int_{|\eta'|=1}  \left[
 \frac{x_n}{\sqrt{a^{nn}}} \left(x'\cdot \eta'
 +i x_n \sqrt{\sum_{j,k=1}^{n-1} \frac{a^{jk}\eta_j \eta_k}{a^{nn}}}
 \right)^{1-n}\right] ds_{\eta'}.
\end{eqnarray*}
 (\ref{3..3}) shows that $u(x)$ is uniquely determined by
the data of $h$ on its support set. (\ref{3..2}) can be rewritten as
\begin{eqnarray} \label{3..4} \left\{\begin{array}{ll}
\big(a^{nn}\big)^2 \frac{\partial^4 u}{\partial x_n^4}
+2\,a^{nn} \frac{\partial^2}{\partial x_n^2}\left(\sum_{j,k=1}^{n-1} a^{jk}
\frac{\partial^2}{\partial x_j \partial x_k} \right)u
\\
 \,\, \quad \quad \quad \; \quad \; +\left(\sum_{j,k=1}^{n-1} a^{jk}
 \frac{\partial^2}{\partial x_j \partial x_k}
 \right)^2 u=0\quad \, &\mbox{in}\;\; {\Bbb R}^n_+, \\
   u=0 \;\; \quad \quad \quad & \mbox{on}\;\; \partial {\Bbb R}^n_+,\\
  \sqrt{a^{nn}} \frac{\partial u}{\partial x_n}=h \quad \quad
 &  \mbox{on}\;\; \partial {\Bbb R}^n_+.\end{array}\right.\end{eqnarray}
 Taking the Fourier transform to
 (\ref{3..4}) with respect to $x_1, \cdots, x_{n-1}$,
 we have
\begin{eqnarray} \label{3.2} \left\{\begin{array}{ll}
\big(a^{nn}\big)^2 \frac{\partial^4 \hat{u}}{\partial x_n^4} - 2 a^{nn}
 \left(\sum_{j,k=1}^{n-1} a^{jk} \eta_j\eta_k\right)\frac{\partial^2 \hat{u}}{\partial x_n^2}
\\
 \quad \quad \quad \quad \quad +\left(\sum_{j,k=1}^{n-1} a^{jk} \eta_j\eta_k\right)^2
  \hat{u}=0 \quad \; & \mbox{in}\;\; {\Bbb R}^n_+,\\
  \hat{u}=0 \quad \quad \;\quad\quad \;  &\mbox{on}\;\; \partial {\Bbb R}^n_+, \\
 \sqrt{a^{nn}} \frac{\partial \hat{u}}{\partial x_n}=\hat{h}(\eta')
 \,\quad \; &\mbox{on}\;\; \partial {\Bbb R}^n_+.\end{array}\right.\end{eqnarray}
 We denote $|\xi'|:=\sqrt{\frac{\sum_{j,k=1}^{n-1} a^{jk} \eta_j\eta_k}{a^{nn}}}$.
 Then, the general solution of (\ref{3.2}) has the form:
\begin{eqnarray} \hat{u}(\eta',x_n)= C_1 e^{|\xi'|x_n}+C_2 e^{-|\xi'|x_n}+C_3 x_n
e^{|\xi'|x_n}+C_4 x_n e^{-|\xi'|x_n},\end{eqnarray} where $C_1$,
$C_2$, $C_3$, $C_4$ are arbitrary functions in $\eta'$. From the boundary
conditions of (\ref{3.2}), it follows that
\begin{eqnarray*}
\hat{u}(\eta', x_n)
&=&C_2(-e^{|\xi'|x_n}+e^{-|\xi'|x_n}+2 x_n  |\xi'|e^{-|\xi'|x_n})\\
  && +C_3(x_ne^{|\xi'|x_n}-x_n
e^{-|\xi'|x_n})+\frac{\hat{h}(\eta')}{\sqrt{a^{nn}}}
x_ne^{-|\xi'|x_n}.\end{eqnarray*}
 Therefore
\begin{eqnarray*} u(x)= \frac{1}{(2\pi)^{n-1}}\int_{{\Bbb R}^n} e^{i\langle x', \eta'\rangle}\left[
C_2 (-e^{|\xi'|x_n}+ e^{-|\xi'|x_n}+2x_n |\xi'|e^{-|\xi'|x_n})\right.\\
\quad \left.+C_3(x_ne^{|\xi'|x_n} -x_n
e^{-|\xi'|x_n})+\frac{\hat{h}(\eta')}{\sqrt{a^{nn}}} x_n
e^{-|\xi'|x_n}\right]d\eta',\end{eqnarray*} from which we have
\begin{eqnarray*} \frac{\partial^2 u}{\partial x_j\partial x_k}&=&
\frac{1}{(2\pi)^{n-1}}\int_{{\Bbb R}^{n-1}} e^{i\langle x', \eta'\rangle} (-\eta_j\eta_k)
\left[C_2 (-e^{|\xi'|x_n}+e^{-|\xi'|x_n}  +2 x_n |\xi'|e^{-|\xi'|x_n}) \right.\\
  && \left.   +C_3(x_ne^{|\xi'|x_n} -x_n
e^{-|\xi'|x_n})+\frac{\hat{h}(\eta')}{\sqrt{a^{nn}}}\,x_ne^{-|\xi'|x_n}\right]d\eta',\end{eqnarray*}
\begin{eqnarray*} \frac{\partial^2 u}{\partial x_n^2}
&=& \frac{1}{(2\pi)^{n-1}}\int_{{\Bbb R}^{n-1}}  e^{i\langle x', \eta'\rangle}
 \left[ C_2 \left(-|\xi'|^2e^{|\xi'|x_n}+|\xi'|^2e^{-|\xi'|x_n}-4|\xi'|^2
e^{-|\xi'|x_n} \right.\right.  \\
&& \;\, \left. \left. +2x_n |\xi'|^3 e^{-|\xi'|x_n}\right)
 +C_3\left(2|\xi'|e^{|\xi'|x_n}+x_n
|\xi'|^2e^{|\xi'|x_n}+2|\xi'|e^{-|\xi'|x_n}
-x_n|\xi'|^2e^{-|\xi'|x_n}\right)\right. \\
&&\quad \; \left. -2\frac{\hat{h}(\eta')}{\sqrt{a^{nn}}}|\xi'|e^{-|\xi'|x_n}
+ \frac{\hat{h}(\eta')}{\sqrt{a^{nn}}} \, x_n |\xi'|^2
e^{-|\xi'|x_n}\right]d\eta'.\end{eqnarray*} Then
\begin{eqnarray*}
&&\left( \sum_{j,k=1}^{n-1} a^{jk} \frac{\partial^2 u}{\partial x_j
\partial x_k} + a^{nn} \frac{\partial^2 u}{\partial
x_n^2}\right)\bigg|_{x_n=0}=
 \left(a^{nn} \frac{\partial^2 u}{\partial
x_n^2}\right)\bigg|_{x_n=0}\\
&& \quad \quad \quad =\frac{1}{(2\pi)^{n-1}}\int_{{\Bbb R}^{n-1}}   a^{nn}
 \left[
-4C_2 |\xi'|^2 + 4C_3|\xi'|-2 |\xi'| \frac{\hat{h}(\eta')}{\sqrt{a^{nn}}}
\right]e^{i\langle x',
\eta'\rangle}\,d\eta'.\end{eqnarray*}
In order to take a bounded solution of
the equation (\ref{3.2}), we may let $C_2=C_3=0$. Hence
\begin{eqnarray*}  &&\left( \sum_{j,k=1}^{n-1}
a^{jk} \frac{\partial^2 u}{\partial x_j \partial x_k} + a^{nn}
\frac{\partial^2 u}{\partial x_n^2}\right)\bigg|_{x_n=0} =
\frac{1}{(2\pi)^{n-1}}\int_{{\Bbb R}^{n-1}} a^{nn}
 \left[ -2|\xi'| \frac{\hat{h}(\eta')}{\sqrt{a^{nn}}}
 \right]  e^{i\langle x', \eta'\rangle}\, d\eta'\quad \quad \quad \\
&&  \quad \quad =\frac{1}{(2\pi)^{n-1}}\int_{{\Bbb R}^{n-1}}  a^{nn}
 \left[ -2\left(\sqrt{\sum_{j,k=1}^{n-1}
  \frac{a^{jk}\eta_j\eta_k}{a^{nn}}}\right)e^{i\langle x',
  \eta'\rangle}\,\frac{\hat{h}(\eta')}{\sqrt{a^{nn}}} \right]
  d\eta'\, \\
 &&\quad \quad =\frac{1}{(2\pi)^{n-1}}\int_{{\Bbb R}^{n-1}}
 \left( -2\sqrt{\sum_{j,k=1}^{n-1} a^{jk}\eta_j\eta_k}\right)
  e^{i\langle x', \eta'\rangle}\, \hat{h}(\eta') d\eta',\qquad\qquad \qquad\end{eqnarray*}
i.e.,
\begin{eqnarray}  F_0h = \frac{1}{(2\pi)^{n-1}}\int_{{\Bbb R}^{n-1}}
 \left(2\sqrt{\sum_{j,k=1}^{n-1}
a^{jk}\eta_j\eta_k}\right) e^{i\langle x', \eta'\rangle}\, \hat{h}(\eta')\,
  d\eta'.\end{eqnarray}
This  shows that the principal symbol of the pseudodifferential operator $F_0$ is
$2\sqrt{\sum_{j,k=1}^{n-1} a^{jk} \eta_j\eta_k}$. $\quad \qquad \square$

\vskip 0.23 true cm

 \noindent {\bf Theorem 3.2} \ \ {\it Let $(\mathcal{M},g)$ be an $n$-dimensional
 $C^\infty$ Riemannian manifold, and let
 $\Omega$ be a bounded domain with $C^\infty$ boundary.
 Assume that the pseudodifferential
 operator $F$ is defined as before.  Then for any coordinate chart
 $\kappa: \partial \Omega \supset U\to U^\kappa \subset {\Bbb R}^{n-1}$
there is a pseudodifferential operator
 $\Lambda\in \Psi^1(U^\kappa)$ such that for every $h\in H^{1/2} (\partial \Omega)$ we have
  \begin{eqnarray*} \kappa_* (F h)-\Lambda(\kappa_* h)\in C^\infty (U^\kappa),\end{eqnarray*}
where $\kappa_*$ is the linear tangent mapping of $\kappa$.
$\Lambda$ is elliptic and has a coordinate invariant positively
$1$-homogeneous principal symbol
 $p_0\in C^\infty (T^* U^\kappa\setminus 0)=C^\infty (U^\kappa \times ({\Bbb R}^{n-1}\setminus 0))$
given by \begin{eqnarray*}  p_0(x', \eta') =2
\sqrt{\sum_{j,k=1}^{n-1} g^{jk} \eta_j\eta_k},\quad \quad \forall
\,\,  (x', \eta')\in T^* U^\kappa\setminus 0.\end{eqnarray*}}

 \noindent  {\it Proof.} \ \  It is well-known (see, for example,
\cite{Mol}) that there is a $T>0$ and a neighborhood $G\subset \mathcal{M}$ of the
 boundary $\partial \Omega$  together with
a diffeomorphism $\psi : G \to \partial \Omega \times [0,T)$ such
that

 i)  \ \  $\,\,\psi(q)=(q,0)$ for every $q\in \partial \Omega$,

 ii) \ \   The unique
geodesic normal to $\partial \Omega$ (with the unit-speed
$\sqrt{g^{nn}}$ with respect to $g$) starting in $q\in
\partial \Omega$ is given by \begin{eqnarray*}  [0,T)\to \Omega,
\;\quad  x_n\to \psi^{-1} (q,x_n).\end{eqnarray*}
 Moreover, $\psi$ is unique with i) and
 ii) and has the following additional properties:
 Let $\kappa: \partial \Omega \supseteq U\to U^\kappa \subset {\Bbb
 R}^{n-1}$  be any coordinate chart on $\partial \Omega$
 and ${\tilde \kappa}: \mathcal{M} \supseteq {\tilde U} \to U^\kappa
 \times [0,T)$ be its extension via $\psi$.

  1)  \ \  The metric $g$ has on $U^\kappa \times [0,T)$ the form
  \begin{eqnarray} \label {3+8} {\tilde \kappa}_* g = \sum_{j,k=1}^{n-1}\big( g_{jk}
  \, dx_j \otimes dx_k\big)
  + g_{nn}dx_n\otimes dx_n, \end{eqnarray}

  2) \ \   The Laplace-Beltrami operator $\Delta_g$ can in $U^\kappa \times [0, T)$
   be written as
\begin{eqnarray*}
\Delta_g^{{\tilde \kappa}}= g^{nn} \,\frac{\partial^2}{\partial
x_n^2} +\frac{1}{\sqrt{|g|}}\, \frac{\partial (\sqrt{|g|}\, g^{nn})}{\partial x_n} \,
\frac{\partial}{\partial x_n} +\sum_{j,k=1}^{n-1}
\frac{1}{\sqrt{|g|}}\
   \frac{\partial}{\partial x_j} \left( \sqrt{|g|}\,g^{jk}
   \frac{\partial}{\partial x_k}\right).
\end{eqnarray*}
 ${\tilde \kappa}$ is called a boundary
normal coordinate chart and its coordinates $x_1, \cdots, x_{n-1},
x_n$ boundary normal coordinates. $G$ is said to be a tubular
neighborhood of $\partial \Omega$.
Therefore, for given $\epsilon>0$ and every  $q\in
\partial \Omega$,
we let
  $(G_{q,\epsilon},  {\tilde \kappa})$
is a boundary normal coordinates chart (Note that $
{\tilde\kappa}(G_{q,\epsilon})=U^\kappa_{q,\epsilon}\times [0,T)$) such that
$\mbox{diam}(G_{q, \epsilon})<\epsilon$. Then there is a partition
of unity subordinate to the open cover $\{G_{q, \epsilon} \cap
\partial \Omega\big| q\in \partial \Omega\}$, i.e., a collection of real-valued $C^\infty$
functions $\phi_i$ on $\partial \Omega$ satisfying the following
conditions:

(a) \  \  The supports of the $\phi_i$ are compact and locally
finite;

(b) \ \  The support of $\phi_i$ is completely contained in
$G_\alpha$ for some $\alpha$;

(c) \ \  The $\phi_i$  sum to one at each point of $\partial
\Omega$:
 \begin{eqnarray*}
\sum_i \phi_i  (x)=1.\end{eqnarray*} Since $h=\sum_i h\phi_i$, we
may assume that the support set of $h$ is contained in some small
neighborhood $G_{q,\epsilon}\cap \partial \Omega$  for some $q\in
\partial \Omega$.  (Let us point out that $\big({\tilde\kappa}(G_{q, \epsilon})\big)\cap \{x_n=0\}
=U^\kappa_{q,\epsilon}\subset \partial {\Bbb
R}^n_+$). It is clear that we can always
choose a fine cover $\{G_{q,\epsilon} \cap \partial \Omega\}$ of
$\partial \Omega$,
 so that in addition to (\ref{3+8}) we have
\begin{eqnarray} \label {302} |\big(g^{jk} (x',x_n)-g^{jk} (0)\big)
\eta_j \eta_k | \le {\epsilon} \sum_{j=1}^{n} \eta_j^2, \quad \;\;
((x',x_n)\in {\tilde \kappa}(G_{q,\epsilon}))
\end{eqnarray} for any given ${\epsilon}>0$, all real $\eta_1,
\cdots, \eta_n$ and all ${\tilde \kappa} (G_{q,\epsilon})$.
  The finer cover does not influence the fact, which is obviously true in the original
 cover, that (\ref{302}) holds.

 Let $u$ be the solution of (\ref{3.1}).
For any fixed $\epsilon_0>0$, in the local coordinates, the $u$ satisfies
 \begin{eqnarray}  \label{3/10} \quad \quad \; \qquad  \quad \left\{\begin{array}{ll}  \left[ \sum_{j,k=1}^{n-1}
\frac{1}{\sqrt{|g|}}\,\frac{\partial}{\partial x_j}\left( \sqrt{|g|}g^{jk} \,\frac{\partial}{\partial x_k}\right)
+ g^{nn} \frac{\partial^2}{\partial x_n^2} \right.\\
   \left. \qquad \qquad \qquad  \qquad +\frac{1}{\sqrt{|g|}} \,\frac{
\partial (\sqrt{|g|} g^{nn})}{\partial x_n} \;\frac{\partial}{\partial x_n} \right]^2 u=0
 \;\;  & \mbox{in}\;\; {\tilde \kappa} (G_{q,\epsilon_0}), \\
u=0 \quad  &\mbox{on}\;\;
 {\tilde
\kappa} (G_{q,\epsilon_0})\cap \partial {\Bbb R}^n_+, \\
 \sqrt{g^{nn}(x)}\,\frac{\partial u}{\partial x_n}=h \;\;
   & \mbox{on}\;\; {\tilde
\kappa} (G_{q,\epsilon_0})\cap \partial {\Bbb R}^n_+. \end{array}\right.\end{eqnarray}
  By the regularity of elliptic equations, we get $u\in C^{4,\alpha} ({\tilde\kappa}(G_{q,\epsilon}))$.
  We define the  operator $F_{q, \epsilon_0}:\, C^\infty_0
(({\tilde\kappa}(G_{q,\epsilon_0}))\cap \partial {\Bbb R}^n_+) \to
 C^\infty (({\tilde\kappa}(G_{q,\epsilon_0}))\cap \partial {\Bbb R}^n_+
 )$ by $F_{q,\epsilon_0} h:=(-\Delta u)\big|_{({\tilde\kappa}(G_{q,\epsilon_0}))
\cap  \partial {\Bbb R}^n_+}$ for any $h\in
C^\infty_0(({\tilde\kappa} (G_{q,\epsilon_0}))\cap \partial {\Bbb R}^n_+
 )$.
 It is easy to check that $F_{q,\epsilon_0}$ is a pseudodifferential
 operator on the open set
$({\tilde\kappa}(G_{q, \epsilon_0})) \cap \partial {\Bbb R}^n_+$. We
denote its principal symbol as $p_0^{(\epsilon_0)} (x', \eta')$.

Next, noticing that (see, for example,  Theorem 1.4.4 of \cite{Jo})
\begin{eqnarray*}\frac{\partial g^{jk}}{\partial x_l}(0)=0 \quad \; \mbox{for all}\;\; 1\le j,k,l\le n,\end{eqnarray*}
 we get
\begin{eqnarray*} \Delta_{g(0)}^2 = \left( \sum_{j,k=1}^{n-1}
g^{jk}(0) \frac{\partial^2 }{\partial x_j \partial x_k} + g^{nn}(0)
\frac{\partial^2 }{\partial x_n^2}\right)^2, \end{eqnarray*}
 where $0$ is the boundary normal coordinates of $q\in \partial \Omega$.
Let $v$ be the solution of the problem
\begin{eqnarray}  \label{3/10}\left\{\begin{array}{ll}  \left( \sum_{j,k=1}^{n-1}
g^{jk}(0) \frac{\partial^2 }{\partial x_j \partial x_k} + g^{nn}(0)
\frac{\partial^2 }{\partial x_n^2}\right)^2 v=0
\quad \;  & \mbox{in}\;\; {\Bbb R}^n_+, \\
v=0 \quad \quad \quad  \quad &\mbox{on}\;\; \partial {\Bbb R}^n_+, \\
 \sqrt{g^{nn}(0)}\frac{\partial v}{\partial x_n}=h \;\quad\;  & \mbox{on}\;\;
 \partial {\Bbb R}^n_+.\end{array}\right.\end{eqnarray}
Since \begin{eqnarray*}
0= \Delta^2_{g(0)} v -\Delta_g^2 u= \Delta_{g(0)}^2(v-u) +(\Delta^2_{g(0)}- \Delta^2_g)u \quad
\; \mbox{in}\;\; {\tilde \kappa} (G_{q,\epsilon}), \end{eqnarray*}
 we get that $v-u$ satisfies
\begin{eqnarray*} \left\{ \begin{array}{ll}
 \Delta_{g(0)}^2(v-u) =(\Delta_g^2- \Delta_{g(0)}^2)u \quad
\; & \mbox{in}\;\; {\tilde \kappa}(G_{q,\epsilon}), \\
v-u=0 \quad \;\; &\mbox{on}\;\;   ({\tilde
\kappa} (G_{q,\epsilon}))\cap \partial {\Bbb R}^n_+ \\
\sqrt{g^{nn}(0)}\, \frac{\partial (v-u)}{\partial x_n} =
(\sqrt{g^{nn}(x)} - \sqrt{g^{nn}(0)}) \frac{\partial u}{\partial x_n}
 \quad \; & \mbox{on}\;\;   ({\tilde
\kappa}(G_{q,\epsilon}))\cap \partial {\Bbb R}^n_+. \end{array} \right. \end{eqnarray*}
It follows from $L^p$-estimates of elliptic equations (see, for example, Theorem 15.3 of \cite{ADN}) that
 \begin{eqnarray*} & \|v-u\|_{W^{4, p}({\tilde \kappa}(G_{q,\epsilon/2}))} \le   C_1
\left(\|(\Delta_g^2 -\Delta^2_{g(0)})u\|_{L^p({\tilde \kappa}(G_{q,\epsilon}))}
 \right.\qquad \qquad \qquad \quad \qquad \quad \qquad \quad \quad \quad \\   & \qquad \quad  \left. + \bigg\| \frac{\sqrt{g^{nn}(x)}- \sqrt{g^{nn}(0)}}{\sqrt{g^{nn}(0)}}\,
  \frac{\partial u}{\partial x_n}\bigg\|_{W^{3-\frac{1}{p}, p}({\tilde \kappa}(G_{q,\epsilon})\cap
  \partial {\Bbb R}^n_+ )}+
  \| v-u\|_{L^p({\tilde \kappa}(G_{q,\epsilon}))}\right), \end{eqnarray*}
where the constant $C_1$ is independent of $v-u$,
$W^{l, p}({\tilde \kappa}(G_{q,\epsilon}))$ is the Sobolev space with $p>n$, and $0<\epsilon<\epsilon_0$.
  From
  \begin{eqnarray*}  \lim_{\epsilon\to 0} \, \|v-u\|_{L^p ({\tilde \kappa}
(G_{q, \epsilon}))} =0, \quad \;   \lim_{\epsilon\to 0} \, \|(\Delta_g^2 -\Delta^2_{g(0)})u\|_{L^p ({\tilde \kappa}
(G_{q, \epsilon}))} =0 \end{eqnarray*}
and
\begin{eqnarray*} \lim_{\epsilon\to 0} \, \bigg\|\frac{\sqrt{g^{nn} (x)} - \sqrt{g^{nn}(0)}}{g^{nn}(0)}
\frac{\partial u}{\partial x_n}\bigg\|_{W^{3-\frac{1}{p}, p}  ({\tilde \kappa}
(G_{q, \epsilon})\cap \partial {\Bbb R}^n_+)}=0,  \end{eqnarray*}
 we obtain
  \begin{eqnarray*} \lim_{\epsilon\to 0} \|v-u\|_{W^{4, p} ({\tilde \kappa}(G_{q,\epsilon/2}))}=0. \end{eqnarray*}
Combining this and the Sobolev imbedding theorem, we find
 \begin{eqnarray*} \|v-u\|_{C^{3,\alpha} ({\tilde \kappa}(G_{q, \epsilon/2}))}=0.\end{eqnarray*}
 Furthermore, applying the Schauder estimates (see, for example,
  the proof of Theorem 7.2 in \cite{ADN}), we have that for any $\epsilon<\epsilon_0$,
 \begin{eqnarray*}  && \|v-u\|_{C^{4,\alpha}({\tilde \kappa} (G_{q,\epsilon/4}))}
\le C_2\left(\|(\Delta^2_g -\Delta^2_{g(0)})u\|_{C^\alpha ({\tilde \kappa} (G_{q, \epsilon/2}))}\right.\quad \quad
\\  && \qquad \quad \quad \quad  \left.  + \bigg\| \frac{\sqrt{g^{nn}(x)}- \sqrt{g^{nn}(0)}}{\sqrt{g^{nn}(0)}}\,
  \frac{\partial u}{\partial x_n}\bigg\|_{C^{3, \alpha}({\tilde \kappa}(G_{q,\epsilon/2}))\cap \partial {\Bbb R}^n_+}
+\|v-u\|_{C^\alpha ({\tilde \kappa} (G_{q, \epsilon/2}))}\right),\end{eqnarray*}
 where the constant $C_2$ is independent of $v-u$.
 By \begin{eqnarray*} \lim_{\epsilon\to 0} \|(\Delta^2_g -\Delta^2_{g(0)})u\|_{C^\alpha ({\tilde
\kappa} (G_{q, \epsilon/2}))}  =0\end{eqnarray*}
 and \begin{eqnarray*} \lim_{\epsilon\to 0} \, \bigg\|\frac{\sqrt{g^{nn} (x)} - \sqrt{g^{nn}(0)}}{g^{nn}(0)}
\frac{\partial u}{\partial x_n}\bigg\|_{C^{3,\alpha}  ({\tilde \kappa}
(G_{q, \epsilon})\cap \partial {\Bbb R}^n_+)}=0,  \end{eqnarray*}
 we have  \begin{eqnarray*}\lim_{\epsilon\to 0} \|v-u\|_{C^{4,\alpha}
({\tilde \kappa}(G_{q,\epsilon/4}))}=0,\end{eqnarray*}
 which implies
 \begin{eqnarray*}\lim_{\epsilon \to 0} \|\Delta_g u-\Delta_{g(0)}v \|_{C^{2,\alpha} ({\tilde \kappa} (G_{q,\epsilon/4}))}
=0. \end{eqnarray*}
 Applying the trace theorem,
we obtain  \begin{eqnarray*}
 \lim_{\epsilon\to 0} \left((\Delta_{g} u)\big|_{C(\tilde{\kappa}(G_{q,\epsilon/4})\cap
 \partial {\Bbb R}^n_+)} -
   (\Delta_{g(0)} v)\big|_{C(\tilde{\kappa}(G_{q,\epsilon/4})\cap
 \partial {\Bbb R}^n_+)}\right)=0,\end{eqnarray*}
    i.e.,
    \begin{eqnarray}  \label{0-0}
  (\Delta_{g} u)\big|_{\{x'=0\big|(x', 0)\in {\Bbb R}^{n-1}\}}
  = (\Delta_{g(0)}v )\big|_{\{x'=0\big|(x', 0)\in {\Bbb R}^{n-1}\}}.\end{eqnarray}

    It follows from Lemma 3.1 that $F_0$ has the principal symbol $2\sqrt{\sum_{j,k=1}^{n-1}
 g^{jk}(0)\eta_j\eta_k}$, where $F_0: \, C_0^\infty({\tilde \kappa}(G_{q,\epsilon_0})\cap \partial {\Bbb R}^n_+)\to
 C^\infty({\tilde \kappa}(G_{q,\epsilon_0})\cap \partial {\Bbb R}^n_+)$ defined by
 $F_0 h = (\Delta v)\big|_{{\tilde \kappa}(G_{q,\epsilon_0})\cap \partial {\Bbb R}^n_+}$,
 and $v$ is the solution of (\ref{3/10}).
  From (\ref{0-0}), we get that  the principal symbol of $F_{q, \epsilon_0}$ is also
$p_0^{(\epsilon_0)} (0,\xi')=2\sqrt{\sum_{j,k=1}^{n-1} g^{jk}(0)
\eta_j\eta_k}$.
  Since $q$ is an arbitrary point at $\partial \Omega$, it follows that the principal symbol of $F$ is
  $2\sqrt{\sum_{j,k=1}^{n-1} g^{jk}(x') \eta_j\eta_k}$ in the local boundary
normal coordinate system $(x')$ of  $\partial \Omega$. $\quad  \qquad \square$

\vskip 0.38 true cm

 3.2. \  We define another pseudodifferential operator $\Theta: H^{3/2} (\partial \Omega) \to
 H^{-3/2} (\partial \Omega)$  as follows:
  Let $\phi \in H^{3/2} (\partial \Omega)$, and  let $u\in H^2(\Omega)$  be the solution of
 \begin{eqnarray} \label{3.,1}  \left\{ \begin{array} {ll}
\triangle_g^2 u=0 \,\quad &\mbox{in}\;\; \Omega,\\
  u=\phi \quad \;\; \quad  & \mbox{on}\;\; \partial \Omega, \\
 \frac{\partial u}{\partial \nu}=0 \quad  & \mbox{on}\;\; \partial
 \Omega,\end{array}\right. \end{eqnarray}
  we set $\Theta h:= \frac{\partial (\Delta_g u)}{\partial \nu}\big|_{\partial \Omega}$.
  By Green's formula, we have
  $$\langle \Theta h, h\rangle=\int_{\partial \Omega}
(\Theta h)h\, ds =\int_{\partial \Omega} u\left(\frac{\partial (\Delta_g u)}{\partial \nu} \right)ds =
\int_\Omega |\triangle_g u|^2 dx\ge 0,\;\; \,
\mbox{for any}\;\; h\in H^{1/2} (\partial \Omega), $$
 which implies that $\Theta$ is a non-negative, self-adjoint pseudodifferential
  operator from $H^{3/2}(\partial \Omega)$ to $H^{-3/2} (\partial \Omega)$.

\vskip 0.25 true cm

\noindent {\bf Lemma 3.3} \ \ {\it  Let $A$ be a positive definition, real-valued constant
 matrix as in (\ref{2--1}).  Assume that
\begin{eqnarray*}  \Theta_0: C^\infty_0({\Bbb R}^{n-1}) \to C^\infty ({\Bbb R}^{n-1})\end{eqnarray*}
  defined by the following problem:  Let $\phi\in C^\infty_0 ({\Bbb R}^{n-1})$ and
let $u\in C^\infty({\Bbb R}^n_+)$ be the  solution of
\begin{eqnarray} \label{3,.2} \left\{\begin{array}{ll} \left( \sum_{j,k=1}^{n-1}
a^{jk} \frac{\partial^2}{\partial x_j \partial x_k} + a^{nn}
\frac{\partial^2 }{\partial x_n^2}\right)^2 u=0
\quad \;\;  & \mbox{in}\;\; {\Bbb R}^n_+, \\
 u=\phi \quad \quad  & \mbox{on}\;\; \partial {\Bbb R}^n_+, \\
 \frac{\partial u}{\partial x_n}=0 \;\quad
  & \mbox{on}\;\; \partial {\Bbb R}^n_+.\end{array}\right.\end{eqnarray}
 we set $\Theta_0\phi:= \sqrt{a^{nn}}\,\frac{\partial (\Delta_g u)}{\partial x_n}\big|_{\partial {\Bbb R}^n}$.
   Then  the principal symbol of $\Theta_0$ is
 \begin{eqnarray*}  p_0(x', \eta') =2\,
\sqrt{\left(\sum_{j,k=1}^{n-1} a^{jk} \eta_j\eta_k\right)^3},\quad \quad \forall
\,\,  (x', \eta')\in {\Bbb R}^{n-1}\times ({\Bbb R}^{n-1}\setminus 0).\end{eqnarray*}}

\vskip 0.2 true cm

\noindent  {\it Proof.} \ \  Similar to Lemma 3.2,
 it follows from Lemma 2.5 that
 \begin{eqnarray} \label {3,.3} u(x', x_n)= \int_{{\Bbb R}^{n-1}} K_1(x'-y', x_n) \phi(y') dy',\end{eqnarray}
where
 \begin{eqnarray*}
 K_1(x', x_n) &=&
  (-1)^{n-1}\frac{(n-2)!}{(2\pi i)^{n-1}}\int_{|\eta'|=1}
 \left[\left(x'\cdot\eta' +i x_n \sqrt{\sum_{j,k=1}^{n-1}
 \frac{a^{jk}\eta_j\eta_k}{a^{nn}}}\right)^{1-n} \right.\\
 && \left. \,\, + (n-1)i x_n \sqrt{\sum_{j,k=1}^{n-1} \frac{a^{jk}\eta_j\eta_k }{a^{nn}}}
  \left(x'\cdot \eta' +i x_n
 \sqrt{\sum_{j,k=1}^{n-1}
 \frac{a^{jk}\eta_j\eta_k}{a^{nn}}}\right)^{-n}\right]ds_{\eta'}.
\end{eqnarray*}
  Taking the Fourier transform for (\ref{3,.2}) with respect to $x_1, \cdots, x_{n-1}$,
 we have
\begin{eqnarray} \label{3,.4} \left\{\begin{array}{ll}
\big(a^{nn}\big)^2 \frac{\partial^4 \hat{u}}{\partial x_n^4} - 2 a^{nn}
 \left(\sum_{j,k=1}^{n-1} a^{jk} \eta_j\eta_k\right)\frac{\partial^2 \hat{u}}{\partial x_n^2}
\\
 \quad \quad \quad \quad \quad +\left(\sum_{j,k=1}^{n-1} a^{jk} \eta_j\eta_k\right)^2
  \hat{u}=0 \quad \; & \mbox{in}\;\; {\Bbb R}^n_+,\\
  \hat{u}=\hat{\phi}(\eta') \quad \quad \;\quad\quad \quad \;\quad &  \mbox{on}\;\; \partial {\Bbb R}^n_+, \\
  \frac{\partial \hat{u}}{\partial x_n}=0
 \,\quad \; & \mbox{on}\;\; \partial {\Bbb R}^n_+.\end{array}\right.\end{eqnarray}
  Then,  the general solution of (\ref{3,.4}) has the form:
\begin{eqnarray} \hat{u}(\eta',x_n)= C_1 e^{|\xi'|x_n}+C_2 e^{-|\xi'|x_n}+C_3 x_n
e^{|\xi'|x_n}+C_4 x_n e^{-|\xi'|x_n},\end{eqnarray} where
 $|\xi'|:=\sqrt{\frac{\sum_{j,k=1}^{n-1} a^{jk} \eta_j\eta_k}{a^{nn}}}$, and
$C_1$, $C_2$, $C_3$, $C_4$ are arbitrary functions in $\eta'$.
In order to obtain a bounded solution of (\ref{3,.4}), we put $C_1=C_3=0$,
so that we find by the boundary
conditions of (\ref{3,.4}) that
 $C_2=\hat{\phi}(\eta')$, $C_4= \hat{\phi}(\eta')|\xi'|$.
 That is, \begin{eqnarray*}
\hat{u}(\eta', x_n)
= \hat{\phi} (\eta')e^{-|\xi'| x_n} \left( 1
 + x_n |\xi'| \right).\end{eqnarray*}
  Thus
\begin{eqnarray*} u(x)&=& \frac{1}{(2\pi)^{n-1}}\int_{{\Bbb R}^{n-1}}  e^{i\langle x', \eta'\rangle}
\hat{\phi} (\eta')  e^{-|\xi'| x_n} \left(
  1+  x_n |\xi'|\right) d\eta'.\end{eqnarray*}
  Since \begin{eqnarray*}  \frac{\partial^3 u}{\partial x_j \partial x_k\partial x_n}
  &=&\frac{1}{(2\pi)^{n-1}}\int_{{\Bbb R}^{n-1}} x_n e^{i\langle x', \eta'\rangle }\eta_j \eta_k \hat{\phi}
  (\eta') |\xi'|^2 e^{-|\xi'| x_n} d\eta'\\
        \frac{\partial^3 u}{\partial x_n^3}
  &=&\frac{1}{(2\pi)^{n-1}}\int_{{\Bbb R}^{n-1}} e^{i\langle x', \eta'\rangle }\hat{\phi}
  (\eta')e^{-|\xi'|x_n}\left( 2|\xi'|^3
    -x_n |\xi'|^4  \right)  d\eta',
  \end{eqnarray*}
  it follows that
  \begin{eqnarray*} \Theta_0\phi &=&  \left[\sqrt{a^{nn}}
  \frac{\partial}{\partial x_n} \left(\sum_{j,k=1}^{n-1}
a^{jk} \frac{\partial^2 u}{\partial x_j \partial x_k} + a^{nn}
\frac{\partial^2 u}{\partial x_n^2}\right)\right]\bigg|_{x_n=0}\\
 &=&
 \frac{1}{(2\pi)^{n-1}}\int_{{\Bbb R}^{n-1}}
  2 \left(\sqrt{a^{nn}}|\xi'|\right)^3  e^{i\langle x', \eta'\rangle} \hat{\phi} (\eta')\,
  d \eta' \\
    &=& \frac{1}{(2\pi)^{n-1}}\int_{{\Bbb R}^{n-1}}
 2 \, \sqrt{\left(\sum_{j,k=1}^{n-1}
 a^{jk}\eta_j\eta_k\right)^3}\, \, e^{i\langle x', \eta'\rangle} \hat{\phi}(\eta')\,
  d\eta'.\end{eqnarray*}
This  shows that the principal symbol of the pseudodifferential operator
$\Theta_0$ on $\partial {\Bbb R}^n_+$ is
$2 \,\sqrt{\left(\sum_{j,k=1}^{n-1} a^{jk} \eta_j\eta_k\right)^3}$. $\quad \qquad \square$

  \vskip 0.3 true cm

 \noindent {\bf Theorem  3.4} \ \ {\it Let $(\mathcal{M},g)$ be an $n$-dimensional
 $C^\infty$ Riemannian manifold, and let
 $\Omega$ be a bounded domain with $C^\infty$ boundary.
 Assume that the pseudodifferential
 operator $\Theta$ is defined as before.  Then for any coordinate chart
 $\kappa: \partial \Omega \supset U\to U^\kappa \subset {\Bbb R}^{n-1}$
there is a pseudodifferential operator
 $\Upsilon\in \Psi^1(U^\kappa)$ such that for every $\phi\in H^{3/2} (\partial \Omega)$ we have
  \begin{eqnarray*} \kappa_* (\Theta \phi)-\Upsilon(\kappa_* \phi)\in C^\infty (U^\kappa).\end{eqnarray*}
$\Upsilon$ is elliptic and has a coordinate invariant positively
$3$-homogeneous principal symbol
 $p_0\in C^\infty (T^* U^\kappa\setminus 0)=C^\infty (U^\kappa \times ({\Bbb R}^{n-1}\setminus 0))$
given by \begin{eqnarray*}  p_0(x', \eta') =2
\left(\sum_{j,k=1}^{n-1} g^{jk} \eta_j\eta_k\right)^{3/2},\quad \quad \forall
\,\,  (x', \eta')\in T^* U^\kappa\setminus 0.\end{eqnarray*}}

\vskip 0.28 true cm

\noindent  {\it Proof.} \ \ The proof is similar to Theorem 3.2.

\vskip 1.49 true cm

\section{Proofs of main results}

\vskip 0.45 true cm

\noindent  {\bf Proof  of Theorem 1.1.} \ \
   Let $F: \, H^{1/2}(\partial \Omega)\to H^{-1/2}
 (\partial \Omega)$ is defined as in \S 3.1.
  It follows from the discussion in \S 3.1 that $F$ is a self-adjoint, elliptic, pseudodifferential
 operator on $H^{1/2}(\partial \Omega)$  whose principal symbol
 is $2\sqrt{\sum_{j,k=1}^{n-1} g^{jk}(x') \eta_j \eta_k}$,
 where $x'$ is the local boundary normal coordinate on $\partial \Omega$.
We define the operator $Z_\epsilon$ by $$Z_\epsilon f(x')
=\left(\frac{1}{\varrho(x')+\epsilon}\right) f(x') \quad \mbox{for
all}\;\; f\in H^{1/2} (\partial \Omega) \;\; \mbox{and}\;\; x'\in \partial \Omega,$$ where $\epsilon>0$ is a sufficiently
small constant.
 Applying Lemma 2.4, we obtain that the operator $Q_\epsilon =Z_\epsilon\circ F
  \,:\, H^{1/2}(\partial \Omega)\to H^{-1/2} (\partial \Omega)$
  defined by $Q_\epsilon h=
\left(\frac{1}{(\varrho(x')+\epsilon)}\, (-\triangle_g
u)\right)\big|_{\partial \Omega}$ is a pseudodifferential operator
 with the principal symbol
$\frac{2\,\sqrt{\sum_{j,k=1}^{n-1} g^{jk} (x')
\eta_j\eta_k}}{\varrho(x')+\epsilon}$,
 where $u$ is the solution of (\ref{3.1}).
  It is easily seen that the operator $Q_\epsilon$ has the same  eigenvalues
$\lambda_k(\epsilon)$ and corresponding normalized  eigenfunctions $\frac{\partial u_k}{\partial \nu}$ on $\partial \Omega$
    as the
 following biharmonic Steklov eigenvalue problem:
 \begin{eqnarray*} \left\{ \begin{array}{ll}  \triangle_g^2 u_k=0
 \;\;\quad & \mbox{in}\;\; \;  \Omega,\\
   u_k=0 \; \quad \quad \; \; &\mbox{on}\;\; \partial \Omega,\\
   \triangle_g u_k
 +(\lambda_k(\epsilon)) (\varrho(x')+\epsilon) \frac{\partial u_k}{\partial \nu}=0\quad\; & \mbox{on}\;\; \partial
 \Omega.\end{array} \right. \end{eqnarray*}
 Let $E_\tau$ be the spectral resolution of $Q_\epsilon$, and let $e(x', y', \tau)$ be the kernel of
$E_\tau$ (here $e(x',y',
\tau)= \sum_{\lambda_k\le \tau}
 \sqrt{(\varrho(x')+\epsilon)(\varrho(y')+\epsilon)}\,\frac{\partial u_k(x')}{\partial
\nu} \, \frac{\partial u_k(y')}{\partial \nu}$).
 It follows from Lemma 2.6 (see, Theorem 1.1 of \cite{Ho}, or \cite{Ho3}) that
     \begin{eqnarray} \tau^{-(n-1)} e(x',x',\tau)
-(2\pi)^{-(n-1)} \int_{B_{x'}} d\xi^* =O(\tau^{-1}) \,\quad
\,\mbox{as}\;\; \tau\to +\infty,\end{eqnarray} where $B_{x'} = \{
\eta'\in T^*_{x'}(\partial \Omega)\big| p_0(x',\eta')<1\}$, $p_0(x',
\eta')$ denotes the principal symbol of $Q_\epsilon$.
 By  $A_\epsilon(\tau)=\int_{\partial \Omega} e(x',x', \tau) dx'$ and
 $p_0(x',\eta')=\frac{2\sqrt{\sum_{j,k=1}^{n-1} g^{jk}(x')
  \eta_j\eta_k}}{\varrho(x')+\epsilon}$, we obtain
 that
\begin{eqnarray*}  A_\epsilon(\tau)= \frac{1}{(2\pi)^{n-1}} \left(\int_{\partial \Omega} dx'
  \int_{2(\varrho(x')+\epsilon)^{-1}\big(\sqrt{\sum_{j,k=1}^{n-1} g^{jk}(x')\, \eta_j\eta_k}\big)<1}
d\xi^*\right) \tau^{n-1}\\
\quad \quad \quad \quad  +O(\tau^{n-2}) \quad \;\; \quad
  \,\mbox{as}\;\; \tau\to +\infty.\qquad\qquad \quad \qquad  \end{eqnarray*}
For any fixed local boundary normal coordinate $x'\in {\tilde \kappa}(\partial \Omega)$,
since $(n-1)\times (n-1)$ matrix $g'=(g^{jk}(x'))$ is positive definite,
there exists an $(n-1)\times (n-1)$ matrix $C(x')=(c_{jk}(x'))$ such that
${}^tC(x') g'(x') C(x')= (\delta_{jk})$, where $\delta_{jk}$ is  the Kronecker
delta. Note that $d\xi^*=\sqrt{|g'(x')|}\,d\zeta_1\cdots
d\zeta_{n-1}$, in each fiber of $T^*(\partial \Omega)$,
which is a vector space of dimension $n-1$. With the change of variables
    $\eta_j=\sum_{k=1}^{n-1} c_{jk}(x') \zeta_k$,
we obtain
\begin{eqnarray*}  && \int_{2(\varrho(x')+\epsilon)^{-1}
\big(\sqrt{\sum_{j,k=1}^{n-1} g^{jk}(x') \eta_j\eta_k}\big)<1}
d\xi^* \qquad\quad  \\
&&\quad  \quad \quad = \int_{\{(\zeta_1, \cdots, \zeta_{n-1})\in
{\Bbb R}^{n-1}\big| \sqrt{\zeta_1^2 +\cdots
  +\zeta_{n-1}^2} <\frac{\varrho(x')+\epsilon}{2}\}}
 |\mbox{det}\, C(x')| \, \sqrt{|g'(x')|} \,d\zeta_1\cdots d\zeta_{n-1})\\
 && \quad \quad\quad =   \int_{\{(\zeta_1, \cdots, \zeta_{n-1})\in {\Bbb R}^{n-1}\big| \sqrt{\zeta_1^2
  +\cdots +\zeta_{n-1}^2} <\frac{\varrho(x')+\epsilon}{2}\}}
 d\zeta_1\cdots d\zeta_{n-1}\\
 && \quad\quad \quad =   \omega_{n-1}\left(\frac{\varrho(x') +\epsilon}{2}\right)^{n-1}, \end{eqnarray*}
here we have used the fact that $|\mbox{det}\, C(x')| \sqrt{|g(x')|}=1$, and where
$\omega_{n-1}$ is the volume of the unit ball of ${\Bbb R}^{n-1}$.
Therefore
 \begin{eqnarray*}  A_\epsilon(\tau)= \frac{1}{(2\pi)^{n-1}}
\omega_{n-1}\,\tau^{n-1} \int_{\partial \Omega} \left(\frac{\varrho(x') +\epsilon}{2}\right)^{n-1} dx'
 +O(\tau^{n-2})  \quad \;
\mbox{as}\;\; \tau\to +\infty.\end{eqnarray*} Letting
$\epsilon\to 0$, we obtain
\begin{eqnarray*} A(\tau)=\frac{1}{(2\pi)^{n-1}}\,\omega_{n-1}\tau^{n-1}
 \int_{\partial \Omega} \left(\frac{\varrho(x')}{2}\right)^{n-1} dx' +O(\tau^{n-2})
 \quad \;
\mbox{as}\;\; \tau\to +\infty,\end{eqnarray*}
 that is,
\begin{eqnarray*} A(\tau)= \frac{\omega_{n-1}\tau^{n-1}}{(4\pi)^{n-1}}
 \int_{\partial \Omega} (\varrho (s))^{n-1} ds +O(\tau^{n-2})
 \quad \; \mbox{as}\;\; \tau\to +\infty.\end{eqnarray*}
 \qed

\vskip 0.5 true cm

\noindent  {\bf Proof  of Theorem 1.2.} \ \
   Let $R_\epsilon: \, H^{3/2}(\partial \Omega)\to H^{-3/2}
 (\partial \Omega)$ be defined as follows:
 For any $\phi\in H^{3/2}(\partial \Omega)$, we put $R_\epsilon\phi
 =\left(\frac{1}{(\varrho+\epsilon)^3}\,\frac{\partial (\triangle_g v)}{\partial \nu}
 \right)\big|_{\partial \Omega}$, where $v$
  satisfies \begin{eqnarray*}\left\{ \begin{array}{ll} \triangle_g^2 v=0 \;\; &\mbox{in}\;\; \Omega,\\
     u=\phi \;\; & \mbox{on}\;\; \partial
  \Omega,\\
\frac{\partial v}{\partial \nu}=0\;\; & \mbox{on}\;\; \partial \Omega,\end{array}\right.\end{eqnarray*}
 and $\epsilon>0$ is a sufficiently small constant.
  Clearly, $R_\epsilon$ is a  self-adjoint,  elliptic, non-negative pseudodifferential
  operator of order $3$.
  By Lemmas 2.4 and Theorem 3.4, we get that the principal symbol of $R_\epsilon$
  has the form $\frac{2}{(\varrho +\epsilon)^{3}}
 \left(\sum_{j,k=1}^{n-1} g^{jk} \eta_j \eta_k\right)^{3/2}$,
  where $\eta'\in {\Bbb R}^{n-1}$.
    It is easily seen that the operator $R_\epsilon$ has the same eigenvalues $\mu^3_k(\epsilon)$
  and corresponding normalized eigenfunctions $v_k$ on $\partial \Omega$ as the
 following biharmonic Steklov eigenvalue problem:
 \begin{eqnarray*} \left\{ \begin{array}{ll}  \triangle_g^2 v_k=0
 \;\; & \mbox{in}\;\; \Omega,\\
  \frac{\partial v_k}{\partial \nu}=0 \; & \mbox{on}\;\; \partial \Omega,\\
   \frac{\partial (\triangle_g v_k)}{\partial \nu} -\mu_k^3(\epsilon)
   (\varrho+\epsilon)^3 v_k=0\;\; & \mbox{on}\;\; \partial
  \Omega.\end{array} \right. \end{eqnarray*}
It follows from Lemma 2.6 (Theorem 1.1 of \cite{Ho}, also see \cite{Ho3} or \cite{AGMT}) that
\begin{eqnarray} \tau^{-(n-1)/3} e(x',x',\tau)
-(2\pi)^{-(n-1)} \int_{B_{x'}} d\xi^* =O(\tau^{1/3}) \,\quad
\,\mbox{as}\;\; \tau\to +\infty,\end{eqnarray} where $B_{x'} = \{
\eta'\in T^*_{x'}(\partial \Omega)\big| p_0(x',\eta')<1\}$, $p_0(x',
\eta')$ denotes the principal symbol of $F_\epsilon$.
 Since  $B_\epsilon(\tau)=\int_{\partial \Omega} e(x',x', \tau^3) dx'$ and
 $p_0(x',\eta')=\frac{2\sqrt{(\sum_{j,k=1}^{n-1} g^{jk}(x')
  \eta_j\eta_k)^3}}{(\varrho(x')+\epsilon)^3}$, we have
\begin{eqnarray*}  B_\epsilon(\tau)= \frac{1}{(2\pi)^{n-1}} \left(\int_{\partial \Omega} dx'
  \int_{2(\varrho(x')+\epsilon)^{-3}\big(\sqrt{\sum_{j,k=1}^{n-1} g^{jk}(x')\, \eta_j\eta_k}\big)^{3/2}<1}
d\xi^*\right) \tau^{n-1}\\
\quad \quad \quad \quad  +O(\tau^{n-2}) \quad
  \,\mbox{as}\;\; \tau\to +\infty,\qquad\qquad \quad \qquad \quad \end{eqnarray*}
 where $d\xi^*=\sqrt{|g'(x')|}\,d\zeta_1\cdots
d\zeta_{n-1}$. With the change of variables
    $\eta_j=\sum_{k=1}^{n-1} c_{jk}(x') \zeta_k$,
where $(n-1)\times (n-1)$ matrix $C(x')=(c_{jk}(x'))$ satisfies
${}^tC(x') g'(x') C(x')= (\delta_{jk})$, we obtain
\begin{eqnarray*}  && \int_{2(\varrho(x')+\epsilon)^{-3}
\big(\sqrt{\sum_{j,k=1}^{n-1} g^{jk}(x') \eta_j\eta_k}\big)^{3/2}<1}
d\xi^* \qquad\quad  \\
&&\quad  \quad \quad = \int_{\{(\zeta_1, \cdots, \zeta_{n-1})\in
{\Bbb R}^{n-1}\big| \sqrt{\zeta_1^2 +\cdots
  +\zeta_{n-1}^2} <\frac{\varrho(x')+\epsilon}{\sqrt[3]{2}}\}}
 \,d\zeta_1\cdots d\zeta_{n-1}\\
  && \quad\quad \quad =   \omega_{n-1}\left(\frac{\varrho(x') +\epsilon}{\sqrt[3]{2}}\right)^{n-1}, \end{eqnarray*}
which implies  \begin{eqnarray*}  B_\epsilon(\tau)= \frac{1}{(2\pi)^{n-1}}
\omega_{n-1}\,\tau^{n-1} \int_{\partial \Omega} \left(\frac{\varrho(x') +\epsilon}{\sqrt[3]{2}}\right)^{n-1} dx'
 +O(\tau^{n-2})  \quad \;
\mbox{as}\;\; \tau\to +\infty.\end{eqnarray*} Letting
$\epsilon\to 0$, we obtain
\begin{eqnarray*} B(\tau)=\frac{1}{(2\pi)^{n-1}}\,\omega_{n-1}\tau^{n-1}
 \int_{\partial \Omega} \left(\frac{\varrho(x')}{\sqrt[3]{2}}\right)^{n-1} dx' +O(\tau^{n-2})
 \quad \;
\mbox{as}\;\; \tau\to +\infty,\end{eqnarray*}
i.e.,
\begin{eqnarray*} A(\tau)=\frac{\omega_{n-1} \tau^{n-1}}{(\sqrt[3]{16}\,
\pi)^{n-1}}
 \int_{\partial D} \varrho^{n-1}(s) ds +O(\tau^{n-2})
 \quad \; \mbox{as}\;\; \tau\to +\infty.\end{eqnarray*}
 \qed

\vskip 1.39 true cm

\section{The asymptotic formulas are sharp}

\vskip 0.45 true cm

  H\"{o}rmander (see \cite{Ho} or \cite{Sh}) proved that, for a pseudodifferential operator of order $m$
   with principal symbol $p_0(x, \xi)$,
\begin{eqnarray} \bigg|e(x, x, \tau)- (2\pi)^{-n} \int_{p_0(x,\xi)<\tau} d\xi^* \bigg|\le C(1+|\tau|)^{\frac{n-1}{m}}
 \end{eqnarray} uniformly in compact subsets of  $\mathcal{M}$, where $C$ is independent of $x, \tau$.
 Applying this result to our cases, we immediately obtain
\begin{eqnarray} \label{5/2} \bigg|A(\tau)-
     \frac{\omega_{n-1}\tau^{n-1}}{(4\pi)^{n-1}}
 \int_{\partial \Omega} \varrho^{n-1} ds \bigg|\le C (1+|\tau|)^{n-2}, \end{eqnarray}
\begin{eqnarray} \label {5/3} \bigg|B(\tau)-
\frac{\omega_{n-1}\tau^{n-1}}{(\sqrt[3]{16}\pi)^{n-1}}
 \int_{\partial D} \varrho^{n-1} ds \bigg| \le C (1+|\tau|)^{n-2}. \end{eqnarray}

  In this section, we shall show that
  (\ref{5/2}) and (\ref{5/3}) cannot further be improved. More precisely,
  we shall show by two counterexamples (letting $\Omega$ respectively be the unit ball of ${\Bbb R}^n$
  and the unit disk of of ${\Bbb R}^2$,and $\varrho\equiv 1$) that
  the asymptotic formulas (\ref{1-7}) and (\ref{1;-6;}) are the ``best possible''.

  \vskip 0.2 true cm

 First we give some well-known facts concerning
spherical harmonics (See e.g. M\"{u}ller \cite{Mu}). When $\Omega=B$
 we may explicitly determine all the
biharmonic Steklov eigenvalues of (\ref{1-1}). In fact, for each
integer $m\ge 0$, let ${\mathcal{P}}_m({\Bbb R}^n)$ denote the set
of homogeneous polynomials of degree $m$ in $n$ variables, i.e., the
set of functions $u$ of the form
\begin{eqnarray*} u(x)= \sum_{|\alpha|=m} a_\alpha x^\alpha \quad \; \mbox{for}\;\;
x\in {\Bbb R}^n,\end{eqnarray*} with coefficients $a_\alpha\in {\Bbb
C}$. {\it A solid spherical harmonic of degree} $m$
 is an element of the subspace
\begin{eqnarray} {\mathcal{H}}_m({\Bbb R}^n)= \{u\in {\mathcal{P}}_m ({\Bbb R}^n)\big|
\Delta u=0\;\;\mbox{on}\;\; {\Bbb R}^n\}.\end{eqnarray} Let
\begin{eqnarray} N(n, m)=\mbox{dim}\,{\mathcal {H}}_m ({\Bbb R}^n)\quad \,
\mbox{for}\;\; n\ge 1 \;\; \mbox{and}\;\; m\ge 0.\end{eqnarray} Note
that ${\mathcal{P}}_0={\mathcal{H}}_0$ is just the space  of
constant functions, and ${\mathcal{P}}_1={\mathcal{H}}_1$ is just
the space of homogeneous linear functions,
 so
 \begin{eqnarray} N(n, 0)=1 \quad \mbox{and}\;\; N(n, 1)=n\quad \;
 \mbox{for}\;\; n\ge 1.\end{eqnarray}
 It follows from p.$\;$251-252 of \cite{Mc} that
\begin{eqnarray*} N(1, m)=\left\{ \begin{array}{ll} 1 \quad \,
\mbox{if}\;\; m=0\;\; \mbox{or}\;\; 1,\\
0 \quad \, \mbox{if}\;\; m\ge 2,\end{array}\right.\end{eqnarray*}
\begin{eqnarray*} N(2, m)=\left\{ \begin{array}{ll} 1 \quad \,
\mbox{if}\;\; m=0,\\
2 \quad \,  \mbox{if}\;\; m\ge 1,\end{array}\right.\end{eqnarray*}
and
\begin{eqnarray*} N(n, m)=\frac{2m+n-2}{n-2} \left( \begin{array}{ll} m+n-3\\
 \quad n-3\end{array}\right) \quad \;\mbox{for}\;\; n\ge 3\;\; \mbox{and}\;\;
 m\ge 0.\end{eqnarray*}

 The following Lemma was obtained by
 Ferrero, Gazzola and Weth (see, Theorem 1.3 of \cite{FGW}):

\vskip 0.10 true cm

\noindent  {\bf Lemma 5.1.} \ \  {\it If $n\ge 2$ and $\Omega=B$,
then for all $m=0, 1,2,3,\cdots$:

(i) \ \ \  \ the eigenvalues of (\ref{1-1}) are
${\tilde\lambda}_m=n+2m$;

 (ii)  \ \ \  the multiplicity of ${\tilde\lambda}_m$ equals $N(n,m)$;

 (iii) \ \  for all ${\tilde\psi}_m\in {\mathcal{H}}_m ({\Bbb R}^n)$, the
 function ${\tilde\phi}_m(x):= (1-|x|^2) {\tilde\psi}_m(x)$ is an eigenfunction
 corresponding to ${\tilde\lambda}_m$.}

\vskip 0.22 true cm

Now, let $0< \lambda_1\le \lambda_2 \le \cdots \le \lambda_k\le
\cdots$ be all biharmonic Steklov eigenvalues for the ball $B$. From the
above lemma and the formula \begin{eqnarray*} \sum_{j=1}^m \left(\begin{array}{ll}a+j-1\\
\quad \;\; j
\end{array}\right)=\frac{(a+1)(a+2)\cdots
(a+m)}{m!}-1,\end{eqnarray*}
 we get that for $n\ge 2$,
\begin{eqnarray*}   A({\tilde\lambda}_m) &=& \#\{i\big|\lambda_i\le {\tilde \lambda}_m\}
=\sum_{k=0}^m N(n, k) \\ &=& 1+ \sum_{k=1}^m
\frac{(n+2k-2)(n+k-3)!}{(n-2)!\, k!}=1+
\sum_{k=1}^m \left[\frac{(n+k-2)!}{(n-2)!\, k!} + \frac{(n+k-3)!}{(n-2)!(k-1)!}\right]\\
&=&1+\sum_{k=1}^m \left(\begin{array}{ll}n+k-2\\ \;\; \quad \,  k
\end{array}\right) +\frac{(n-2)!}{(n-2)!} +\sum_{k=2}^m
\left(\begin{array}{ll}n+k-3\\ \quad k-1 \end{array}\right) \\
 &=& 1+\sum_{k=1}^m \left(\begin{array}{ll}(n-1)+k-1\\ \quad \quad \quad k \end{array}\right) +1+
\sum_{k=2}^m
\left(\begin{array}{ll}(n-1)+(k-1)-1\\ \quad \quad \quad k-1 \end{array}\right) \\
  &=&2+ \sum_{k=1}^m \left(\begin{array}{ll}(n-1)+k-1\\ \quad \quad \quad k \end{array}\right)+
\sum_{j=1}^m \left(\begin{array}{ll}(n-1)+j-1\\ \quad \quad \quad \;
 j \end{array}\right) -
\left(\begin{array}{ll}(n-1)+m-1\\ \quad \quad \quad  m \end{array}\right) \\
 &=& 2+\left[ \frac{n(n+1)\cdots (n-1+m)}{m!}-1\right] +
\left[\frac{n(n+1)\cdots (n-1+m)}{m!} -1\right] -\left(\begin{array}{ll} n+m-2\\
\quad \;\; m \end{array}\right) \\
 &=& 2 \left[\frac{n(n+1)\cdots (n-1+m)}{m!} \right] -
\left(\begin{array}{ll} n+m-2\\
\,\quad \;\; m \end{array}\right)\\
&=& 2\left(\begin{array}{ll} n+m-1\\  \;\; \quad m\end{array}\right)
 - \left(\begin{array}{ll} n+m-2\\  \; \;\quad m\end{array}\right)
=2\left(\begin{array}{ll} n+m-1\\  \, \quad n-1\end{array}\right) -
\left(\begin{array}{ll} n+m-2\\  \quad  \, n-2\end{array}\right)
.\end{eqnarray*} By applying the formula
\begin{eqnarray*} \left(\begin{array}{ll} p+1\\   \quad
r\end{array}\right) - \left(\begin{array}{ll} \quad  p\\
 r-1\end{array}\right) = \left(\begin{array}{ll}  p\\
r\end{array}\right)\end{eqnarray*} and $m=\frac{{\tilde
\lambda}_m}{2}-\frac{n}{2}$, we obtain
\begin{eqnarray*} A({\tilde\lambda}_m) &=&
\left(\begin{array}{ll} n+m-1\\  \, \quad n-1\end{array}\right) +
\left(\begin{array}{ll} n+m-2\\  \quad\, n-1\end{array}\right)
\\ &=& \frac{(2m+n-1)(m+n-2)(m+n-3)\cdots (m+1)}{(n-1)!}\\
&=& \frac{1}{2^{n-2}\; (n-1)!} ({\tilde \lambda}_m-1)
 ({\tilde \lambda}_m+n-4)({\tilde \lambda}_m+n-6)
\cdots ({\tilde \lambda}_m-n+2). \end{eqnarray*}
 In view of $\frac{\sqrt{\pi} \, \Gamma
(2z)}{2^{2z-1}}=\Gamma(z)\Gamma(z+\frac{1}{2})$,
we  get \begin{eqnarray*} \frac{1}{2^{n-2}\, (n-1)!}&=& \frac{2^{n} \pi^{n-1}}{ (4\pi)^{n-1} \, (n-1)!}=
 \frac{n \, \pi^{n-\frac{1}{2}}}{(4\pi)^{n-1} \,\frac{\pi^{\frac{1}{2}}
\Gamma(n+1)}{2^n}}\\ &=& \frac{n \,
\pi^{n-\frac{1}{2}}}{(4\pi)^{n-1} \Gamma
(\frac{n}{2}+\frac{1}{2})\Gamma(\frac{n}{2}+1)}
=\frac{1}{(4\pi)^{n-1}}\cdot
\frac{\pi^{\frac{n-1}{2}}}{\Gamma(\frac{n-1}{2}+1)}\cdot
\frac{n\pi^{\frac{n}{2}}}{\Gamma(\frac{n}{2}+1)}\\ &=&
\frac{\omega_{n-1}}{(4\pi)^{n-1}}\cdot
n\omega_n=\frac{\omega_{n-1}}{(4\pi)^{n-1}} (\mbox{vol} (\partial
B)).\end{eqnarray*}
 Hence \begin{eqnarray*}  A({\tilde\lambda}_m)&=&
\frac{\omega_{n-1}}{(4\pi)^{n-1}} (\mbox{vol} (\partial B)) ({\tilde
\lambda}_m-1)  ({\tilde \lambda}_m+n-4)({\tilde \lambda}_m+n-6)
\cdots ({\tilde \lambda}_m-n+2)\\
 &=& \frac{\omega_{n-1}}{(4\pi)^{n-1}} (\mbox{vol} (\partial B))
 \left[ {\tilde \lambda}_m^{(n-1)} +(1-n) {\tilde \lambda}_m^{n-2} +
\cdots - (n-4) (n-6) \cdots (-n+2)\right].\end{eqnarray*}
 Since $1-n\ne 0$, this shows that the formula (\ref{1-7}) is sharp.

\vskip 0.25 true cm

 Similar to the above $A({\tilde \lambda}_m)$, we can also give a counter-example to show that the
 remainder term estimate in the  asymptotic formula (\ref{1;-6;}) is sharp.
 Let $B \subset {\Bbb R}^2$ be the unit disk.  If $m\ge 1$, then the functions
 \begin{eqnarray}  \label{5z1} \psi_{m,1}(r, \theta) =r^m \cos m\theta \quad \, \mbox{and}\quad \,
    \psi_{m,2}(r, \theta) =r^m \sin m\theta  \end{eqnarray}
     form an orthogonal basis of harmonic function in the space ${\mathcal{H}}_m({\Bbb R}^2)$.
     Let us consider the following Neumann boundary value problem:
     \begin{eqnarray}  \label{507} \left\{ \begin{array}{ll}
     \Delta u_{m,j} = \psi_{m, j}  \quad \, & \mbox{in} \;\;  B,\\
     \frac{\partial u_{m,j}}{\partial \nu} =0 \quad \, & \mbox{on}\;\; \partial B,\end{array}\right.
     \end{eqnarray}
     where $m=1,2,3,\cdots; j=1,2$.
     We claim that the above solutions $u_{m,j}$ satisfy
      \begin{eqnarray} \label {55} \left\{ \begin{array}{ll}
     \Delta^2 u_{m,j} = 0  \quad \, & \mbox{in} \;\;  B,\\
     \frac{\partial u_{m,j}}{\partial \nu} =0 \quad \, & \mbox{on}\;\; \partial B,\\
          u_{m,j}=\frac{1}{\mu_{m,j}^3}\, \frac{\partial(\Delta u_{m,j})}{\partial \nu}
          \quad \, & \mbox{on}\;\; \partial B.  \end{array}\right.
     \end{eqnarray}
     where \begin{eqnarray*} \frac{1}{\mu_{m,j}^3} &=& \frac{\int_B |\Delta u_{m,j}|^2 dx}{\int_{\partial B}
     \left(\frac{\partial (\Delta u_{m,j})}{\partial \nu}\right)^2 ds}
     = \frac{\int_B |\psi_{m,j}|^2 dx}{\int_{\partial B}
     \left(\frac{\partial \psi_{m,j}}{\partial \nu}\right)^2 ds}=\frac{1}{2m^2(m+1)}.\end{eqnarray*}
  In fact, for the Dirichlet problem
   \begin{eqnarray*}  \left\{ \begin{array}{ll} \Delta u_{m,j} =\psi_{m,j}(r, \theta) \,  &\mbox{in}\;\; B,\\
    u_{m,j} =\frac{-1}{2m(m+1)}  \psi_{m,j} (1, \theta) \, &\mbox{on}\;\; \partial B,
 \end{array} \right. \quad \; m=1,2,3,\cdots; \, j=1,2,\end{eqnarray*}
 from the formula of the solution to the Dirichlet
 boundary value problem for the Poisson equation in the unit disk, we have
\begin{eqnarray*} u_{m,j} (\rho, \alpha) =-\frac{1}{4\pi}
\int_0^1 r\,dr \int_0^{2\pi} \ln \frac{1+\rho^2 r^2 -2r \rho \cos (\alpha-\theta)}
  {r^2 +\rho^2 -2r\rho \cos (\alpha -\theta)} \psi_{m,j} (r, \theta)d\theta  \\
 - \frac{1}{2\pi}  \int_0^{2\pi}  \frac{1-\rho^2}{1+\rho^2 -2 \rho \cos (\alpha-\theta)} \left(
\frac{\psi_{m,j} (1, \theta)}{2m(m+1)} \right)d\theta. \quad \quad\;\; \end{eqnarray*}
Then
\begin{eqnarray*} \frac{\partial u_{m,j}(\rho, \alpha)}{\partial \rho}\big|_{\rho=1}=\frac{1}{2\pi}
\int_0^1  \int_0^{2\pi} \frac{r(1 -r^2)}{1+r^2 -2r\cos (\alpha- \theta)} \psi_{m,j} (r, \theta) dr \, d\theta \\
 - \frac{1}{2\pi}  \left[ \frac{\partial}{\partial \rho} \int_0^{2\pi}  \frac{1-\rho^2}{
 1+\rho^2 -2\rho \cos (\alpha -\theta)} \left(
\frac{\psi_{m,j} (1, \theta)}{2m(m+1)} \right) d\theta\right]\bigg|_{\rho=1}.
\quad \quad \;\; \end{eqnarray*}
 It is well-known that
\begin{eqnarray*} \frac{1-r^2} {1+r^2-2r \cos l(\alpha -\theta)} =1+2\sum_{l=1}^\infty
r^l \cos l(\alpha -\theta).\end{eqnarray*}
Therefore   \begin{eqnarray*} \frac{\partial u_{m,j}(\rho, \alpha)}{\partial \rho}\big|_{\rho=1}&=& \frac{1}{2\pi}
\int_0^1  \int_0^{2\pi}  \left( r+2 \sum_{l=1}^\infty r^{l+1}\cos l (\alpha -\theta)\right) \psi_{m,j} (r,\theta) dr\, d\theta \\
  &&  - \frac{1}{2\pi} \left[ \frac{\partial}{\partial \rho} \int_0^{2\pi}
\left( 1+2 \sum_{l=1}^\infty \rho^l \cos l(\alpha -\theta)\right)
    \frac{\psi_{m,j} (1, \theta)}{2m(m+1)} \, d\theta \right] \bigg|_{\rho=1}
\\  &  = & -\frac{1}{2\pi}
\int_0^1  \int_{\alpha}^{\alpha-2\pi}  \left( r+2 \sum_{l=1}^\infty r^{l+1}\cos l t\right) \psi_{m,j} (r,\alpha-t) dr\, dt
\\   &&   +\frac{1}{2\pi} \left[ \int_\alpha^{\alpha- 2\pi} \left(
   2 \sum_{l=1}^\infty l \rho^{l-1} \cos l t \right)
\frac{\psi_{m,j} (1, \alpha -t)}{2m(m+1)} \, dt \right]\bigg|_{\rho=1}.
\end{eqnarray*}
 By (\ref{5z1}) and a simple calculation, we get $\frac{\partial u_{m,j}}{\partial \nu}\big|_{\rho=1}=0$.
  Combining this and the result
\begin{eqnarray*}  \left\{ \begin{array}{ll} \Delta \psi_{m,j} =0 \,  &\mbox{in}\;\; B,\\
     \frac{\partial \psi_{m,j}}{\partial \nu} = \eta_{m,j} \psi_{m,j} \, &\mbox{on}\;\; \partial B,
 \end{array} \right. \quad \; m=1,2,3,\cdots; \, j=1,2,\end{eqnarray*}
 where  \begin{eqnarray*} \frac{1}{\eta_{m,j}} =\frac{\int_B |\nabla(\psi_{m,j})|^2 dx}{\int_{\partial B}
     \big(\frac{\partial \psi_{m,j}}{\partial \nu}\big)^2 ds}=\frac{1}{m}, \quad \; \, m=1,2,3,\cdots; \; j=1,2,
    \end{eqnarray*}
  we show the desired claim.
    It follows that   \begin{eqnarray*} \mu_{m,j}^3 = 2m^2(m+1), \quad \; m=1,2,3, \cdots; \, j=1,2,\end{eqnarray*}
    Thus, for $n=2$,
\begin{eqnarray*}   B(\mu_{m,j}) &=& \#\{i\big|\mu_i\le \mu_{m,j}\}
=\sum_{k=0}^m N(2, k)=2(m+\frac{1}{2})\\
 &=& \frac{\omega_{1} (\mbox{vol} (\partial B))}{2\pi} \left(m+\frac{1}{2}\right). \end{eqnarray*}
  By \begin{eqnarray*} m\sim \frac{\mu_{m,j}}{\sqrt[3]{2}} \quad \; \mbox{as} \;\; m\to +\infty, \end{eqnarray*}
  we have  \begin{eqnarray*}   B(\mu_{m,j})\sim  \frac{\omega_1 ({\mbox{vol} (\partial B)})}{\sqrt[3]{16}\pi}
  \left(\mu_{m,j} +\frac{1}{2}\right) \quad \; \mbox{as}\;\; m\to +\infty, \,  j=1,2. \end{eqnarray*}
 Since $\frac{1}{2}\ne 0$, this shows that asymptotic formula (\ref{1;-6;})
 cannot be improved on unit disk of ${\Bbb R}^2$.

\vskip 1.68  true cm

\centerline {\bf  Acknowledgments}

\vskip 0.38 true cm

\vskip 0.38 true cm
 I wish to express my sincere gratitude to Professor
 L. Nirenberg, Professor Fang-Hua Lin and Professor J. Shatah for their
  support during my visit at Courant Institute.
  I would like to thank Professor K. C. Chang and  Professor C. Sogge
 for their constructive suggestions to this paper.
  This research was also supported by SRF for ROCS, SEM (No.2004307D01) and NSF of China (No.11171023).

  \vskip 1.65 true cm

\vskip 0.32 true cm

\end{document}